\documentclass{ifacconf}
\usepackage[centertags]{amsmath}
\usepackage{amssymb}
\usepackage{enumitem}
\usepackage{float}
\usepackage{enumerate}
\usepackage{graphicx,subfigure}
\usepackage{nomencl} 
\usepackage{siunitx}
\usepackage{array}
\usepackage[square,sort,comma,numbers]{natbib}
\allowdisplaybreaks
\graphicspath{{./Figs/}}

\newcommand{\scalefont}[1]{\scalebox{1.5}{#1}}

\newcommand{\norm}[1]{\ensuremath{\left\| #1 \right\|}}
\newcommand{\bracket}[1]{\ensuremath{\left[ #1 \right]}}
\newcommand{\braces}[1]{\ensuremath{\left\{ #1 \right\}}}
\newcommand{\parenth}[1]{\ensuremath{\left( #1 \right)}}
\newcommand{\pair}[1]{\ensuremath{\langle #1 \rangle}}

\newcommand{\SO}{\ensuremath{\mathsf{SO(3)}}}
\newcommand{\T}{\ensuremath{\mathsf{T}}}
\renewcommand{\L}{\ensuremath{\mathsf{L}}}
\newcommand{\so}{\ensuremath{\mathfrak{so}(3)}}

\renewcommand{\Re}{\ensuremath{\mathbb{R}}}

\newcommand{\D}{\ensuremath{\mathbf{D}}}

\newcommand{\Sph}{\ensuremath{\mathsf{S}}}
\renewcommand{\S}{\Sph}

\newcommand{\Ad}{\ensuremath{\mathrm{Ad}}}
\newcommand{\ad}{\ensuremath{\mathrm{ad}}}

\newcommand{\G}{\ensuremath{\mathsf{G}}}
\newcommand{\g}{\ensuremath{\mathfrak{g}}}

\begin{document}

\begin{frontmatter}
	
	\title{Geometric Optimal Controls for Flapping Wing UAV on a Lie Group\thanksref{footnoteinfo}} 
	
	\thanks[footnoteinfo]{This research has been supported in part by NSF under the grants NSF CMMI-1761618 and CMMI-1760928.}
	
	\author[First]{Tejaswi K. C.} 
	\author[First]{Taeyoung Lee} 
	
	\address[First]{The George Washington University, Washington DC 20052 (e-mail: kctejaswi999@gmail.com, tylee@gwu.edu).}
	
	\begin{abstract}                
		Inspired by flight characteristics captured from live Monarch butterflies, an optimal control problem is presented while accounting the effects of low-frequency flapping and abdomen undulation. 
       A flapping-wing aerial vehicle is modeled as an articulated rigid body, and its dynamics are developed according to Lagrangian mechanics on an abstract Lie group.
        This provides an elegant, global formulation of the dynamics for flapping-wing aerial vehicles, avoiding complexities and singularities associated with local coordinates.
        This is utilized to identify an optimal periodic motion that minimizes energy variations, and an optimal control is formulated to stabilize the periodic motion.
        Furthermore, the outcome of this paper can be applied to optimal control for any Lagrangian system on a Lie group with a configuration-dependent inertia.
	\end{abstract}
	
	
\end{frontmatter}

\section{Introduction}

Millions of Monarch butterflies migrates from North America to the central Mexico during the fall, exhibiting the longest flight range among the insects (\cite{Gibo1981}).
Their dynamics are distinct from small insects, as the relatively large wings are flapping at a lower frequency, with active undulation of its abdomen. 
It has been suggested that abdomen undulation may reduce power consumption from the dynamic coupling of wing-body motion by \cite{sridhar2019beneficial}.
It is further reported in \cite{dyhr2013flexible} that moths actively modulate their body shape to control flight in response to visual pitch stimuli, and it may contribute to pitch stability.
However, it is challenging to dynamically model such effects to utilize in control system design. 

Flapping wing aerial vehicles are essentially infinite dimensional, nonlinear time-varying systems, where the equations of motion describing displacement and the deformation of a flexible multi-body system are coupled with the Navier-Stokes equations.
Various control system design techniques have been reviewed by~\cite{Shyy}. 
Most of these control systems are based on the common simplified formulation where the nonlinear time-varying flapping dynamics are transformed into linear time-invariant systems by considering small perturbations averaged over the period of flapping (see, for example, ~\cite{Deng2006}).
As such, these approaches are not suitable to analyze the low-frequency flapping dynamics of Monarch butterflies.

Recently, a flapping wing aerial vehicle is modeled as an articulated rigid body by~\cite{sridhar2020geometric}, where four rigid bodies representing two wings, thorax, and abdomen are interconnected via spherical joints, with the assumption of quasi-aerodynamics. 
The resulting dynamics are considered as a Lagrangian system on a Lie group, and an intrinsic form of equations of motion are constructed. 
Compared with developing equations of motion of multi-rigid body systems with local coordinates, such as Euler angles, 
this provides an elegant, global formulation that is free of singularities. 
As such, this is particularly useful to design control systems inspired by Monarch. 
For example, it has been utilized to study the effect of abdomen undulation in energy efficiency by \cite{tejaswi2020effects}.

In this paper, we present an optimal control problem to stabilize a periodic motion representing the hovering flight. 
The flapping motion of both wings are parameterized by several variables characterizing the amplitude and the shape of oscillations, 
which are optimized over the numerical solutions of the aforementioned Lagrangian system. 
Compared with various prior works in the control of flapping wing aerial vehicles, the unique contribution is that we consider the complete dynamics involving the motion of wings, thorax, and abdomen coupled though arbitrary three-dimensional rotations and translations. 
In other words, the dynamics are not simplified by the common assumptions such as the longitudinal motion confined to a two-dimensional space, or the wing flapping decoupled from the body and the abdomen.
These features are particularly useful to grasp the unique dynamic characteristics of Monarch, and to take the advantage of those in control system design. 
In short, we exploit the geometric formulation of Lagrangian mechanics on a Lie group for optimal control of a complex system inspired by Monarch.

\section{Lagrangian Mechanics Formulated on a Lie Group}

Consider  an $n$-dimensional Lie group $\G$.
Let $\g$ be the associated Lie algebra, or the tangent space at the identity, i.e., $\g = \T_e\G$.
Consider a left trivialization of the tangent bundle of the group $\T\G \simeq \G\times \g$, $(g,\dot g)\mapsto (g, \L_{g^{-1}}\dot g)\equiv(g,\xi)$.
More specifically, let $\L:\G\times\G\rightarrow\G$ be the left action defined such that $\L_g h = gh$ for $g,h\in\G$.
Then the left trivialization is a map $(g,\dot g)\mapsto (g, L_{g^{-1}}\dot g)\equiv(g,\xi)$, where $\xi\in\g$.
Further, suppose $\g$ is equipped with an inner product $\pair{\cdot, \cdot}$, which induces an inner product on $\T_g\G$ via left trivialization.
For any $v,w\in\T_g\G$, $\pair{w,v}_{\T_g\G} = \pair{ \T_g \L_{g^{-1}} v, \T_g \L_{g^{-1}} w}_\g$. 
Given the inner product, we identify $\g\simeq \g^*$ and $\T_g \G \simeq \T^*_g \G\simeq G\times \g^*$ via the Riesz representation. Throughout this paper, the pairing is also denoted by the dot product $\cdot$.
The adjoint operator is denoted by $\Ad_g:\g\rightarrow\g$, and the ad operator is denoted by $\ad_\xi:\g\rightarrow\g$. 
See, for example~\cite{MarRat99} for detailed preliminaries. 

We develop Euler--Lagrange equations for an arbitrary Lie group $\G$, which are utilized later for the flapping wing UAV. 

\begin{assum}
	The Lagrangian $L:\G\times \g \rightarrow \Re$ is given by the difference between kinetic and potential energy,
	\begin{equation}\label{eqn:Lg}
		L(g,\xi) = \frac{1}{2} \pair{ \mathbf{J}_g(\xi), \xi } - U(g),
	\end{equation}
	for a configuration-dependent inertia $ \mathbf{J}:\G\times \g\rightarrow \g^* $ and potential $U:\G \rightarrow \Re$.
	
	Here, the inertia is a symmetric, positive-definite tensor dependent on the group. 
    More specifically,
	\begin{gather*}
		\pair{\mathbf{J}_g(\xi), \xi} \geq 0,\\
		\pair{\mathbf{J}_g(\xi), \xi} = 0\; \Leftrightarrow \; \xi=0, \\ 
		\pair{\mathbf{J}_g(\xi_1), \xi_2} = \pair{\mathbf{J}_g(\xi_2), \xi_1},
	\end{gather*}
	for any $g\in\G$ and $\xi,\xi_1,\xi_2\in\g$. 
\end{assum}

\begin{defn}
	The left-trivialized derivative of $ \mathbf{J}_g(\xi) $ with respect to $ g $ is defined as $ \mathbf{K}_g(\xi)(\cdot) :\G\times \g\rightarrow \g^*$ given by
	\begin{equation}
        \mathbf{K}_g(\xi) \chi = \T_e^* \L_g (\D_g \mathbf{J}_g(\xi)) \cdot \chi. \label{eqn:KK}
	\end{equation}
	along the direction $ \chi \in \g $.
	By selecting a basis of $\g$, $\mathbf{K}_g(\xi)$ can be represented by a matrix since it is a linear operator.
\end{defn}

The Euler-Lagrange equations for an arbitrary Lagrangian on a Lie group has been reported, for exampled by \cite{lee2017global}.
Here we present a special case, when the Lagrangian is given by the configuration-dependent kinetic energy and the potential as in \eqref{eqn:Lg}.

\medskip
\begin{thm}\label{thm:EL}
	The forced Euler-Lagrange equations corresponding to the Lagrangian in \eqref{eqn:Lg} are given by,
	\begin{gather}
		\mathbf{J}_g(\dot \xi) + \mathbf{K}_g(\xi)\xi - \ad^*_\xi  \mathbf{J}_g(\xi)  
        - \frac{1}{2}\mathbf{K}^*_g(\xi)\xi \nonumber\\ 
        +\T_e^* \L_g (\D_g U(g)) = f \label{eqn:EL_G}, \\
		\dot g = g \xi,\label{eqn:dot_g}
	\end{gather}
	with $ f: \bracket{t_0, t_f} \rightarrow \g^* $ as the generalized force acting on the system.
\end{thm}
\begin{pf}
	Note that the Lagrangian of the system is expressed as $ L(g, \xi) : \G \times \g \rightarrow \Re $ by utilizing the left-trivialization, $ \xi = g^{-1} \dot{g} $.
	Now, the dynamical relations are obtained using the Lagrange-d'Alembert principle,
	\[
	\delta \int_{t_0}^{t_f} L(g, \xi) dt + \int_{t_0}^{t_f} f(t) \cdot \eta dt = 0
	\]
	where the infinitesimal variation $ \eta = g^{-1} \delta g \in \g $ vanishes at the endpoints.
	Thus the forced Euler-Lagrange equations are (\cite{lee2017global}),
	\begin{equation}\label{eqn:EL_orig}
        \frac{d}{dt} \D_\xi L(g, \xi) - \ad_\xi^* (\D_\xi L(g, \xi)) - \T_e^*\L_g  (\D_g L(g, \xi)) = f.
	\end{equation}
	
	Using the special structure of the Lagrangian from \eqref{eqn:Lg},
	\begin{align*}
		\D_\xi L(g,\xi) \cdot \delta \xi &= \frac{1}{2} \parenth{\pair{\mathbf{J}_g(\xi), \delta \xi} + \pair{\mathbf{J}_g(\delta \xi), \xi}}\\
		&= \pair{\mathbf{J}_g(\xi), \delta \xi}.
    \end{align*}
    Thus, $\D_\xi L(g,\xi) = \mathbf{J}_g(\xi)$.
	We also have,
	\begin{align*}
        \frac{d}{dt} \D_\xi L(g,\xi) & = \mathbf{J}_g(\dot \xi) + \T_e^* \L_g ( \D_g \mathbf{J}_g(\xi)) \cdot \xi \\
		&= \mathbf{J}_g(\dot \xi) + \mathbf{K}_g(\xi)\xi
	\end{align*}
	from the definition in \eqref{eqn:KK}.
	Finally,
	\begin{align*}
		\T^*_e \L_g \cdot \D_g L(g,\xi)\cdot\chi & =
		\frac{1}{2}\pair{\mathbf{K}_g(\xi)\chi , \xi} - \T_e^* \L_g \D_g U(g)\cdot \chi\\
		&= \frac{1}{2}\pair{\chi , \mathbf{K}_g^*(\xi) \xi} - \T_e^* \L_g \D_g U(g)\cdot \chi\\
		& = \parenth{ \frac{1}{2}\mathbf{K}_g^*(\xi)\xi - \T^*_e \L_g \D_g U(g)} \cdot \chi,
	\end{align*}
	where $\g^*$ is identified with $\g$ with the pairing.
	
	Substituting all the above expressions back in \eqref{eqn:EL_orig}, we obtain the Euler-Lagrange equations in \eqref{eqn:EL_G}.
	\hspace*{\fill} \qed
\end{pf}
\medskip



\section{Dynamics of Flapping Wing UAV}\label{sec:dynamics}

In this section, we present a multibody model for an FWUAV after which we derive the corresponding Euler-Lagrange equations.
The three-dimensional special orthogonal group is denoted by $ \SO = \lbrace R \in \mathbb{R}^{3\times3} \mid R^T R = I, \det(R) = 1 \rbrace $, 
and the corresponding Lie algebra is $ \mathfrak{so}(3) = \lbrace A \in \mathbb{R}^{3\times3} \mid A = -A^T \rbrace $.
The \textit{Hat} map $ \wedge : \mathbb{R}^3 \to \mathfrak{s0}(3) $ is defined such that $\hat x y = x\times y$ for any $x,y\in\Re^3$.
And its inverse map is the \textit{vee} map, $ \vee:\so\rightarrow\Re^3 $.
Next, $e_i\in\Re^n$ denotes the $i$-th standard basis of $\Re^n$ for an appropriate dimension $n$, e.g., $e_1=(1,0,\ldots, 0)\in\Re^n$. 
The units are in $\si{kg}$, $\si{m}$, $\si{s}$, and $\si{rad}$, unless specified otherwise.

\subsection{Multibody Model}

\setlength{\unitlength}{0.1\textwidth}
\begin{figure}
    \centerline{
        \footnotesize
        \subfigure[flapping angle,\newline $\phi_R\in[-\pi,\pi)$]{
            \scalebox{0.5}{
                \begin{picture}(3,2.6)(0,0)
                    \put(0,0){\includegraphics[trim={4cm 3cm 4cm 1.5cm},clip,width=0.3\textwidth]{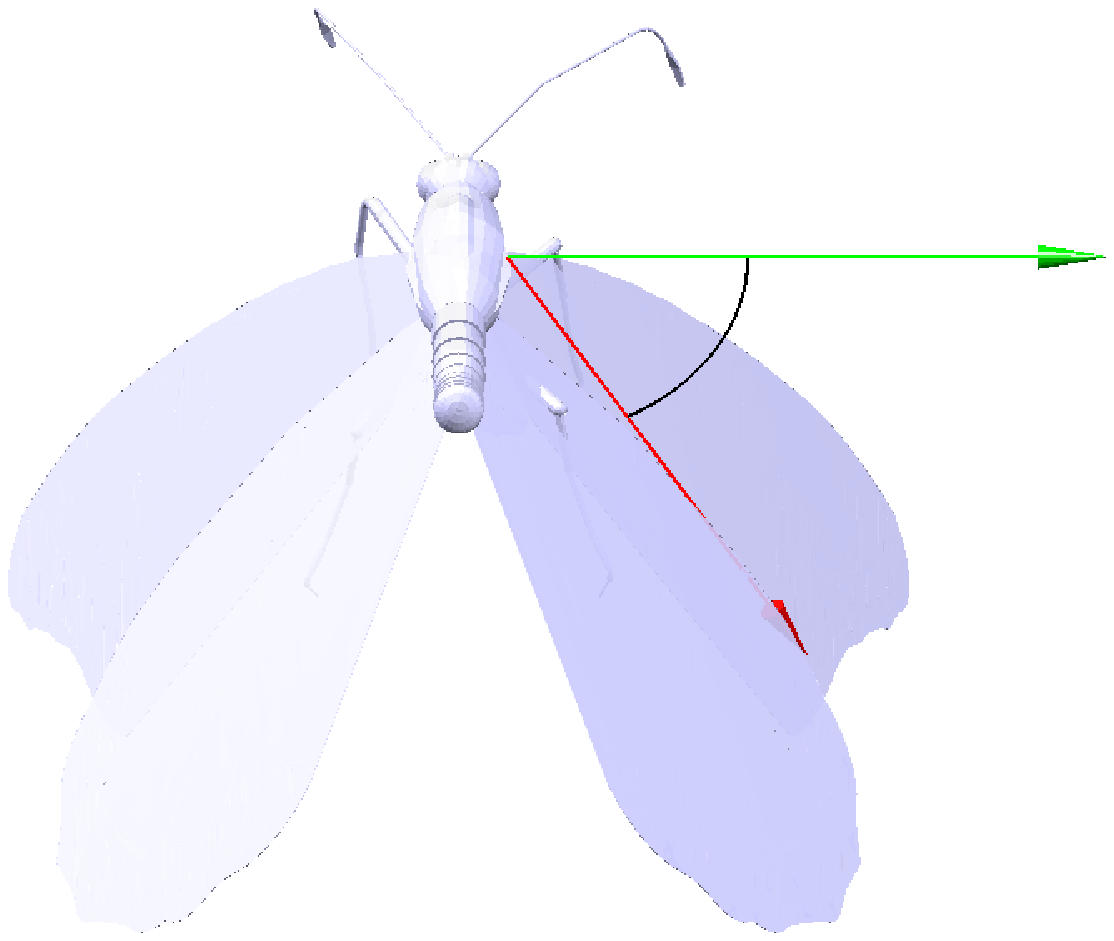}}
                    \put(1.9,1.4){\scalefont{$\phi>0$}}
                    \put(2.75,1.9){\scalefont{$\mathbf{s}_y$}}
                    \put(2.2,0.6){\scalefont{$\mathbf{r}_y$}}
                \end{picture}
            }
        }
        \subfigure[pitch angle,\newline $\theta_R\in[-\pi,\pi)$]{
            \scalebox{0.5}{
                \begin{picture}(3,2.6)(0,0)
                    \put(0,0){\includegraphics[trim={4cm 3cm 2.5cm 1.5cm},clip,width=0.3\textwidth]{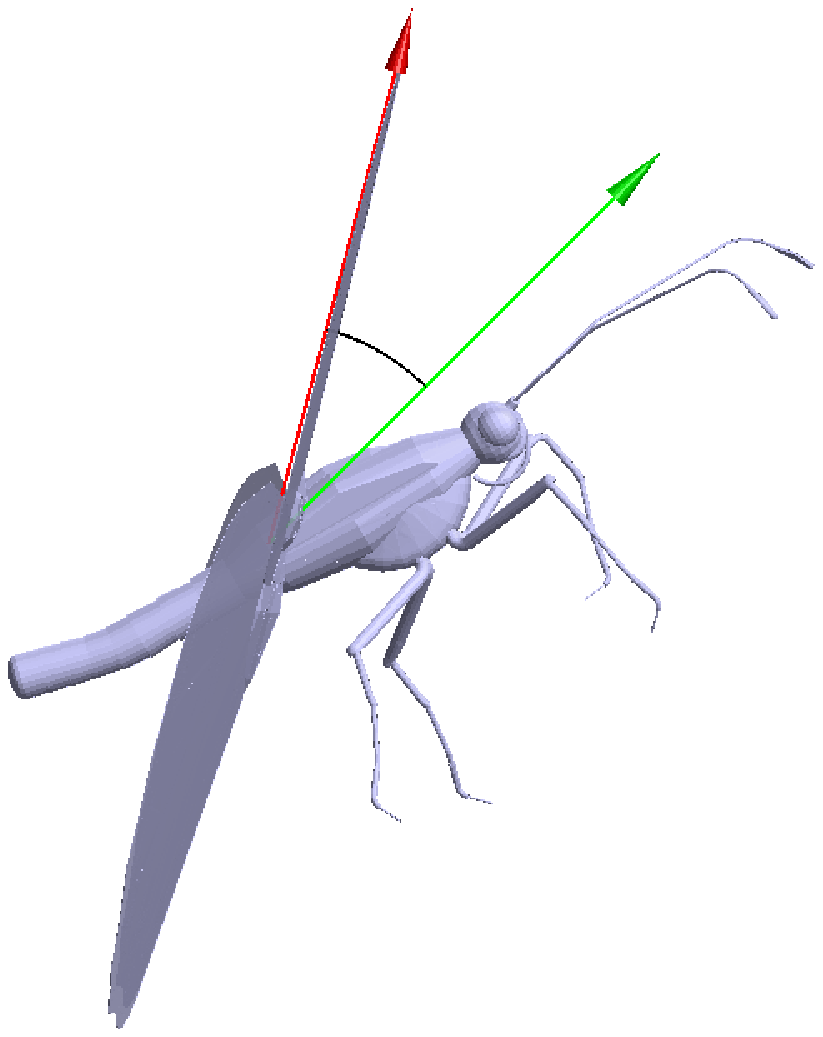}}
                    \put(1.3,1.7){\scalefont{$\theta>0$}}
                    \put(1.1,2.3){\scalefont{$\mathbf{r}_x$}}
                    \put(2.0,1.85){\scalefont{$\mathbf{s}_x$}}
                \end{picture}
            }
        }
        \subfigure[deviation angle,\newline $\psi_R\in[-\pi,\pi)$]{
            \scalebox{0.5}{
                \begin{picture}(3,2.6)(0,0)
                    \put(0,0){\includegraphics[trim={3cm 2cm 3cm 2cm},clip,width=0.3\textwidth]{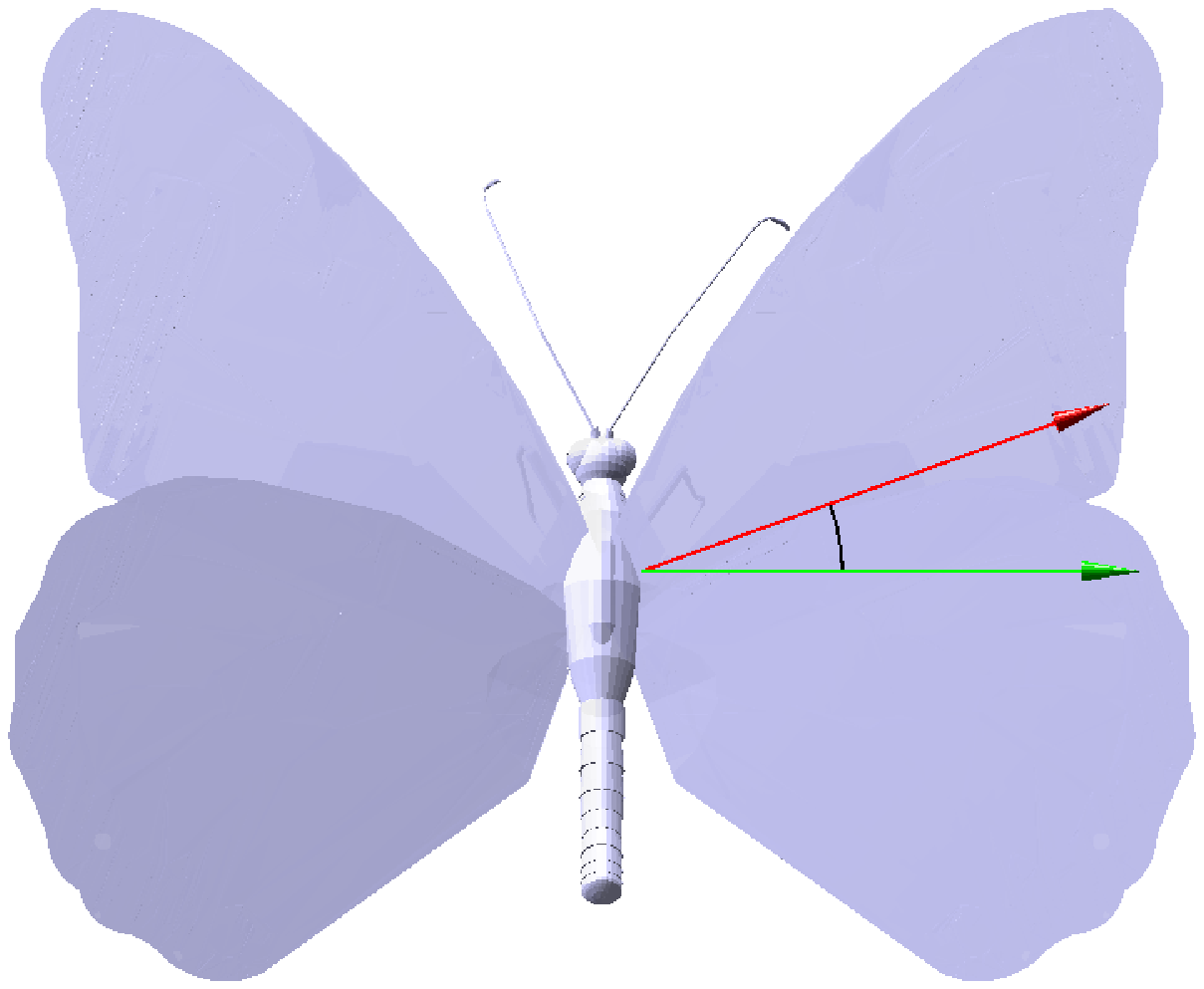}}
                    \put(1.8,1.3){\scalefont{$\psi>0$}}
                    \put(2.8,1.4){\scalefont{$\mathbf{r}_y$}}
                    \put(2.8,0.95){\scalefont{$\mathbf{s}_y$}}
                \end{picture}
            }
        }
    }

	\caption{Euler angles (\cite{sridhar2020geometric}) : positive values are indicated from $ \mathcal{F}_S $ (green) to $ \mathcal{F}_R $ (red)}\label{fig:wing_Euler}
\end{figure}

Let the inertial frame which is compatible to the standard north-east-down (NED) frame be $\mathcal{F}_I=\{\mathbf{i}_x,\mathbf{i}_y,\mathbf{i}_z\}$.
We model the FWUAV as an articulated structure which is composed of multiple rigid bodies listed here: 
\begin{itemize}
	\addtolength\itemsep{3mm}
	\item \textit{Body:} This corresponds to the head and thorax combined into a single rigid body.
	We define $\mathcal{F}_B=\{\mathbf{b}_x,\mathbf{b}_y,\mathbf{b}_z\}$ as the body-fixed frame located at the center of mass of the body.
	This position is denoted by $ x \in \Re^3 $ in $ \mathcal{F}_I $, and the attitude of $\mathcal{F}_B$ is given by $R\in\SO$.
    With $\Omega\in\Re^3$ as the angular velocity of the body resolved in $\mathcal{F}_B$, the attitude evolves as $ \dot R = R \hat \Omega $.

	\item \textit{Right wing:} It is directly attached to the \textit{body}.
	Also, we do not distinguish forewings and hindwings in our model.
	Let $\mathcal{F}_R=\{\mathbf{r}_x,\mathbf{r}_y,\mathbf{r}_z\}$ be the frame fixed to the right wing at its root. 
	It is located at a constant $\mu_R\in\Re^3$ from the origin of $\mathcal{F}_B$.
    Next, we define $\mathcal{F}_S=\{\mathbf{s}_x,\mathbf{s}_y,\mathbf{s}_z\}$ as the stroke frame obtained by translating the origin of $\mathcal{F}_B$ to the center of wing roots, and rotating it about $\mathbf{b}_y$ by a fixed angle $\beta\in[-\pi,\pi)$.
	The attitude of the right wing frame with respect to $\mathcal{F}_S$ is denoted by 1--3--2 Euler angles $(\phi_R(t), \psi_R(t), \theta_R(t))$ (see Figure \ref{fig:wing_Euler}).
    So the attitude of $ \mathcal{F}_R $ relative to $ \mathcal{F}_B $, $Q_R\in\SO$, can be expressed as
	$
	Q_R= \exp(\beta \hat e_2)\exp(\phi_R \hat e_1) \exp(-\psi_R \hat e_3) \exp(\theta_R\hat e_2),
	$
	and its time-derivative will be $ \dot Q_R = Q_R \hat \Omega_R $ for $\Omega_R\in\Re^3$. 

	\item \textit{Left Wing:} The left wing frame $\mathcal{F}_L=\{\mathbf{l}_x, \mathbf{l}_y, \mathbf{l}_z\}$ is defined symmetrically to the right wing, and is located at $\mu_L\in\Re^3$ from the origin of $\mathcal{F}_B$.
    So its attitude is, $ Q_L = \exp(\beta \hat e_2)\exp(-\phi_L \hat e_1) \exp(\psi_L \hat e_3) \exp(\theta_L \hat e_2) ,$
    with the set of Euler-angles $(\phi_L(t), \psi_L(t), \theta_L(t))$, and $ \dot Q_L = Q_L \hat \Omega_L $ for $\Omega_L\in\Re^3$.

	\item \textit{Abdomen:}
	Finally, the abdomen is connected to the body via a spherical joint at which $\mathcal{F}_A=\{\mathbf{a}_x, \mathbf{a}_y, \mathbf{a}_z\}$ is attached.
	It is located at $\mu_A\in\Re^3$ from origin of $\mathcal{F}_B$,
    and its attitude relative to the body is denoted by $Q_A\in\SO$ with $\dot Q_A = Q_A\hat\Omega_A$ for $\Omega_A\in\Re^3$. 
\end{itemize}

\subsection{Lagrangian Mechanics of FWUAV}

The configuration of the presented model is described by  $g=(x, R, Q_R, Q_L, Q_A)$ which belongs to the 15-dimensional Lie group $\G=\Re^3\times \SO^4$. 
The corresponding left trivialized velocity is $\xi = (\dot x, \Omega, \Omega_R, \Omega_L, \Omega_A) $ which is an element of the Lie algebra, $\g = \Re^3 \times \so^4 \simeq \Re^3 \times (\Re^3)^4$. 
In this subsection, we derive the Euler--Lagrange equations in \eqref{eqn:EL_G} for the flapping wing UAV. 

\begin{defn}
	Variables corresponding to the body are denoted by a subscript $ B $, while the wings and the abdomen are denoted by $ \mathcal{B}_i $ with $  i \in \braces{R, L, A} $.
\end{defn}
\medskip
\begin{prop}
	The kinetic energy of the UAV is given by $ \displaystyle T = \frac{1}{2} \xi^T \mathbf{J}_g \xi $ where $ \mathbf{J}_g \in \Re^{15\times 15} $ is the inertia tensor for the UAV given in \eqref{eqn:Jg}.
\end{prop}
\begin{pf}
Firstly, the kinetic energy of the body can be written as,
\begin{align*}
	T_B = \frac{1}{2}m_B \|\dot x\|^2 + \frac{1}{2} \Omega^T J_B\Omega.
\end{align*}
where $m_B\in\Re$ is the mass of the body composed head and thorax, and $J_B\in\Re^{3\times 3}$ is the inertia matrix of the body about $\mathcal{F}_B$.

Next, we need the kinetic energy of the wings and the abdomen which can be obtained in a similar manner.
Consider a mass element $dm$ in $\mathcal{B}_i$, whose location is given by $\nu\in\Re^3$ in $\mathcal{F}_i$.
Thus, its location from the origin of the inertial frame, resolved in $\mathcal{F}_I$ is,
\[
x + R (\mu_i + Q_i \nu) = x + R \mu_i + R Q_i \nu,
\]
and its velocity is 
\[
\dot x + R\hat \Omega (\mu_i + Q_i\nu ) + R Q_i \hat\Omega_i \nu.
\]
Therefore, the kinetic energy corresponding to $ \mathcal{B}_i $ is
\begin{align}\label{eqn:T_i}
	T_i & = \frac{1}{2} \int_{\mathcal{B}_i} \|\dot x + R\hat \Omega (\mu_i + Q_i\nu ) + R Q_i \hat\Omega_i \nu \|^2 dm.
\end{align}
Let $m_i \in \Re$ be the mass of $\mathcal{B}_i$.
Resolved in $ \mathcal{F}_i $, define $ \nu_i\in\Re^3 $ as the location of the mass center of $ \mathcal{B}_i $ and $J_i\in\Re^{3\times 3}$ as the inertia matrix of $\mathcal{B}_i$ about the origin of $\mathcal{F}_i$ :
\begin{align*}
	\nu_i = \frac{1}{m_i} \int_{\mathcal{B}_i} \nu dm, \quad
	J_i = \int_{\mathcal{B}_i} \hat \nu^T \hat \nu  dm.
\end{align*}
Using these expressions, the kinetic energy can be written as
\begin{align*}
	T_i & = \frac{1}{2}
	\begin{bmatrix} \dot x \\ \Omega \\ \Omega_i \end{bmatrix}^T
	\mathbf{J}_i(R,Q_i)
	\begin{bmatrix} \dot x \\ \Omega \\ \Omega_i \end{bmatrix}.
\end{align*}
For instance, $ [\mathbf{J}_i]_{4:6,6:9} $, which can also be denoted as $ \mathbf{J}_{i_{23}} $ in terms of a $ 3 \times 3 $ block structure, can be obtained from \eqref{eqn:T_i} as
\begin{gather*}
	\frac{1}{2} \Omega^T \mathbf{J}_{i_{23}} \Omega_i = \frac{1}{2} \int_{\mathcal{B}_i} \pair{R\hat \Omega (\mu_i + Q_i\nu ), R Q_i \hat\Omega_i \nu} dm \\
	\implies \Omega^T \mathbf{J}_{i_{23}} \Omega_i = \int_{\mathcal{B}_i} -((\hat\mu_i + \widehat{Q_i\nu})\Omega) ^T Q_i \hat\Omega_i \nu dm \\
	= \int_{\mathcal{B}_i} \Omega^T (\hat\mu_i + \widehat{Q_i\nu})^T Q_i \hat\nu \Omega_i dm \\
	\implies \mathbf{J}_{i_{23}} = \int_{\mathcal{B}_i} (\hat\mu_i + \widehat{Q_i\nu})^T Q_i \hat\nu dm = m_i \hat\mu_i^T Q_i \hat\nu_i + Q_i J_i
\end{gather*}
since $ \widehat{Q_i \nu} = Q_i \hat{\nu} Q_i^T $.
Repeating this procedure for all blocks, the configuration-dependent inertia for $\mathcal{B}_i$, i.e., $\mathbf{J}_i(R,Q_i)\in\Re^{9\times 9}$ is given by
{\small
\begin{align}
	\begin{bmatrix}
		m_i I_{3\times 3} & -m_i R(\hat \mu_i + \widehat{Q_i\nu_i}) & - m_i R Q_i \hat \nu_i \medskip \\
		m_i (\hat \mu_i + \widehat{Q_i\nu_i}) R^T & \left(m_i \hat\mu_i^T\hat\mu_i + Q_i J_i Q_i^T +\right. & Q_i J_i + m_i \hat\mu_i^T Q_i \hat\nu_i \\
		& \left. m_i (\hat \mu_i^T \widehat{Q_i\nu_i} + \widehat{Q_i\nu_i}^T \hat\mu_i )\right) & \medskip \\
		m_i\hat\nu_i Q_i^T R^T & J_i Q_i^T + m_i \hat\nu_i^T Q_i^T \hat\mu_i & J_i 
	\end{bmatrix}.
\end{align}
}
The total kinetic energy will be the sum of the individual contributions from the body, the wings and the abdomen
\[
T = \frac{1}{2} \xi^T \mathbf{J}_g \xi = T_B + \sum_{i\in\{R,L,A\}} T_i.
\]
Hence, the symmetric inertia tensor for the complete UAV, $ \mathbf{J}_g \in \Re^{15\times 15} $, can be constructed using the above values as
\begin{align}\label{eqn:Jg}
	\begin{bmatrix}
		\left(m_B I_{3\times 3} + \mathbf{J}_{R_{11}} +\right.
		& \left(\mathbf{J}_{R_{12 }} +\right.
		& \quad \mathbf{J}_{R_{13 }}
		& \mathbf{J}_{L_{13}} 
		& \mathbf{J}_{A_{13}} \\
		\left. \mathbf{J}_{L_{11}} + \mathbf{J}_{A_{11}}\right) & \left.\quad \mathbf{J}_{L_{12}} + \mathbf{J}_{A_{12}}\right) & & & \medskip\\
		\cdot & \quad \left(J_B + \mathbf{J}_{R_{22}} +\right.
		& \quad \mathbf{J}_{R_{23}} 
		& \mathbf{J}_{L_{23}} 
		& \mathbf{J}_{A_{23}} \\
		& \quad \left.\mathbf{J}_{L_{22}} + \mathbf{J}_{A_{22}}\right) & & & \medskip\\
		\cdot & \cdot & \quad\mathbf{J}_{R_{33}}
		& 0 & 0 \\
		\cdot & \cdot & \cdot
		& \mathbf{J}_{L_{33}} & 0 \\
		\cdot & \cdot & \cdot & \cdot & \mathbf{J}_{A_{33}}
	\end{bmatrix},
\end{align}
where $\mathbf{J}_{i_{mn}}, i \in\braces{R,L,A}$ refers to the $m,n$-th $3\times 3$ block of the corresponding matrix, $ \mathbf{J}_i $.
	\hspace*{\fill} \qed
\end{pf}
\medskip

\begin{defn}
	Identify $ \chi = [\delta x, \eta, \eta_R, \eta_L, \eta_A]^T \in \Re^{15} \simeq \g $ as the variation of the configuration, $ g \in \G $.
\end{defn}

Using above expressions of the inertia tensor, we need to evaluate its derivative as defined in \eqref{eqn:KK}.
For example, the first three rows of $ \mathbf{J}_i \cdot [\dot x, \Omega, \Omega_i]^T $ are
\[
m_i \dot x -m_i R(\hat \mu_i + \widehat{Q_i\nu_i}) \Omega - m_i R Q_i \hat \nu_i \Omega_i.
\]
So the corresponding first three rows of $ \mathbf{K}_i(\xi) \chi_i$ along the direction $ \chi_i = [\delta x, \eta, \eta_i]^T \in \Re^9 $ are given by
\begin{gather*}
	[ (\mathbf{K}_i(\xi)) (\chi_i) ]_{1:3} = 
	- m_i R\hat\eta (\hat\mu_i + \widehat{Q_i\nu_i}) \Omega 
	- m_i R( \widehat{Q_i\hat\eta_i\nu_i}) \Omega \\
	- m_i R\hat\eta Q_i \hat\nu_i\Omega_i
	- m_i R Q_i\hat\eta_i \hat\nu_i\Omega_i \\
	= [0] \delta x +  [m_i R ((\hat\mu_i + \widehat{Q_i\nu_i}) \Omega)^\wedge + m_i R (Q_i\hat\nu_i\Omega_i)^\wedge ] \eta \\
	+
	[   - m_i R\hat\Omega Q_i\hat\nu_i
	+ m_i R Q_i \widehat{\hat\nu_i \Omega_i} ] \eta_i.
\end{gather*}
Similarly, repeating this for the remaining rows to construct 
$\mathbf{K}_i(\xi) \in\Re^{9\times 9}$,
{\small
\begin{align}
	\begin{bmatrix}
		0 & m_i R ((\hat\mu_i+\widehat{Q_i\nu_i})\Omega +
		& m_i R( -\hat\Omega Q_i \hat\nu_i + \\
		& Q_i\hat\nu_i\Omega_i)^\wedge  & Q_i \widehat{\hat \nu_i \Omega_i}) \medskip \\
		0 & m_i(\hat\mu_i + \widehat{Q_i\nu_i})\widehat{R^T \dot x}
		& m_i\widehat{R^T\dot x}Q_i\hat\nu_i - Q_i (J_iQ_i^T\Omega)^\wedge + \\
		& & Q_i J_i \widehat{Q_i^T\Omega} - m_i\hat\mu_i\hat\Omega Q_i\hat\nu_i - \\
		& & m_i\widehat{\hat\mu_i\Omega} Q_i \hat\nu_i - Q_i \widehat{J_i\Omega_i} + m_i \hat\mu_i Q_i \widehat{\hat\nu_i \Omega_i} \medskip \\
		0 & m_i \hat\nu_i Q_i^T \widehat{R^T\dot x}
		& m_i \hat\nu_i (Q_i^T R^T \dot x)^\wedge + J_i\widehat{Q_i^T\Omega} - \\
		& & m_i\hat\nu_i (Q^T\hat\mu_i \Omega)^\wedge
	\end{bmatrix}.
\end{align}
}

Thus the derivative of the inertia tensor for the complete UAV from \eqref{eqn:Jg} is expressed as the matrix $\mathbf{K}_g(\xi) \in \Re^{15\times 15}$, 
\begin{align}\label{eqn:Kg}
	\mathbf{K}_g(\xi) = \begin{bmatrix}
		0 & \mathbf{K}_{R_{12}} + \mathbf{K}_{L_{12}} + \mathbf{K}_{A_{12}} & \mathbf{K}_{R_{13}} & \mathbf{K}_{L_{13}} & \mathbf{K}_{A_{13}}\\
		0 & \mathbf{K}_{R_{22}} + \mathbf{K}_{L_{22}} + \mathbf{K}_{A_{22}} & \mathbf{K}_{R_{23}} & \mathbf{K}_{L_{23}} & \mathbf{K}_{A_{23}} \\
		0 & \mathbf{K}_{R_{32}} & \mathbf{K}_{R_{33}} & 0 & 0 \\
		0 & \mathbf{K}_{L_{32}} & 0 & \mathbf{K}_{L_{33}} & 0 \\
		0 & \mathbf{K}_{A_{32}} & 0 & 0 & \mathbf{K}_{A_{33}}
	\end{bmatrix}.
\end{align}

\begin{prop}
	The generalized force due to the gravitational potential energy is,
	\begin{align}\label{eqn:fg}
		\mathbf{f}_g = \begin{bmatrix}
			(m_B+m_R+m_L+ m_A )g e_3 \smallskip\\
			\displaystyle \sum_{i\in\{R,L,A\}} m_i g (\mu_i + Q_i \nu_i)^\wedge {R^T e_3} \smallskip\\
			m_R g \hat\nu_R (Q_R^T R^T e_3) \smallskip\\
			m_L g \hat\nu_L (Q_L^T R^T e_3) \smallskip\\
			m_A g \hat\nu_A (Q_A^T R^T e_3)
		\end{bmatrix}.
	\end{align}
\end{prop}

\begin{pf}
The gravitational potential energy of the body can be written as
\[
U_B = -m_B g e_3^T x,
\]
while that of $\mathcal{B}_i$ is,
\[
U_i = -m_ig e_3^T (x + R\mu_i + RQ_i \nu_i).
\]
So the total potential energy is
\[
U = -m_B g e_3^T x + \sum_{i\in\{R,L,A\}} -m_ig e_3^T (x + R\mu_i + RQ_i \nu_i).
\]

Its negative derivatives, $\mathbf{f}_g\in\Re^{15}$ correspond to the gravitational force and moment given by $ \mathbf{f}_g = - \T^*_e\L_g \D_g U $.
So along the direction $ \chi $,
\begin{align*}
	\delta U &= -m_Bge_3^T \delta x + \\
	& \sum_{i\in\{R,L,A\}} -m_ig e_3^T (\delta x + R\hat{\eta}(\mu_i + Q_i \nu_i) + RQ_i\hat{\eta}_i \nu_i)) \\
	&= -mge_3\cdot \delta x + \sum_{i\in\{R,L,A\}} \bracket{m_ig e_3^T R(\mu_i + Q_i \nu_i)^\wedge} \eta + \\
	& \quad \sum_{i\in\{R,L,A\}} \bracket{m_ig e_3^T RQ_i \hat\nu_i} \eta_i \\
	&= - \mathbf{f}_g \cdot \chi
\end{align*}
from the expression in \eqref{eqn:fg}.
Here, the total mass is denoted by $m\in\Re$,
\begin{align*}
	m = m_B + m_R + m_L + m_A.
\end{align*}
	\hspace*{\fill} \qed
\end{pf}

\begin{prop}
	The contributions from the external aerodynamic forces and control torque are given by
	\begin{align}
		\mathbf{f}_a & = 
		\begin{bmatrix}
			RQ_R F_R  + R Q_L F_L + R Q_A F_A \\
			\hat \mu_R Q_RF_R + \hat\mu_L Q_L F_L + \hat\mu_A Q_A F_A \\
			M_R  \\
			M_L  \\
			M_A   
		\end{bmatrix},\label{eqn:f_a}\\
		\mathbf{f}_\tau & = 
		\begin{bmatrix}
			0 \\
			-\tau_R  -\tau_L  -\tau_A \\
			Q_R^T \tau_R \\
			Q_L^T \tau_L \\
			Q_A^T \tau_A
		\end{bmatrix}\label{eqn:f_tau}
	\end{align}
	where $ F_i $ is the net aerodynamic force and $ M_i $ is the net moment about the wing root or joint connecting the abdomen.
	Moreover, $ \tau_i $ is the control torque exerted at the wing root or abdomen joint resolved in the body-fixed frame.
\end{prop}

\begin{pf}
Consider an infinitesimal aerodynamic force $dF_i(\nu)\in\Re^3$ acting at the position $\nu$ of the wing or the abdomen resolved in the corresponding frame $ \mathcal{B}_i $. 
Thus the net force and moment can be expressed as,
\[
F_i = \int_{\mathcal{B}_i} dF_i(\nu) ,\quad M_i = \int_{\mathcal{B}_i} \nu\times dF_i(\nu).
\]
In the inertial frame, this infinitesimal force is $ R Q_i dF_i(\nu) $ acting at the location $ x+ R\mu_i + R Q_i\nu $. 
So, the corresponding virtual work is
\begin{align*}
	\delta \mathcal{W}_i &= \int_{\mathcal{B}_i} \delta(x + R\mu_i + R Q_i\nu) \cdot R Q_i dF(\nu) \\
	& =  \int_{\mathcal{B}_i} (\delta x + R\hat\eta \mu_i) \cdot R Q_i dF(\nu) + \int_{\mathcal{B}_i} (\hat \eta_i \nu) \cdot dF(\nu) \\
	& = \delta x \cdot R Q_i F_i + (\hat\eta \mu_i) \cdot Q_i F_i + \eta_i \cdot M_i \\
	\implies &\sum_{i\in\{R,L,A\}} \delta \mathcal{W}_i = \mathbf{f}_a \cdot \chi
\end{align*}
Next, the virtual work due to the control torque, which is equal to $ Q_i^T \tau_i $ in the corresponding frame, will be
\begin{align*}
	\delta \mathcal{W}_\tau &=  \sum_{i\in\{R,L,A\}} \eta_i \cdot Q_i^T \tau_i + \eta \cdot -\tau_i = \mathbf{f}_\tau \cdot \chi.
\end{align*}
Here, the second term is the contribution of a reactive torque $ -\tau_i $ exerted on the body.
So the net total of these terms is equal to
\begin{align*}
	\delta\mathcal{W} & =  \delta\mathcal{W}_\tau + \sum_{i\in\{R,L,A\}} \delta\mathcal{W}_i = (\mathbf{f}_{a} + \mathbf{f}_\tau)  \cdot \chi,
\end{align*}
with $\mathbf{f}_a, \mathbf{f}_\tau \in\Re^{15}$ as given in \eqref{eqn:f_a} and \eqref{eqn:f_tau} respectively. 
	\hspace*{\fill} \qed
\end{pf}

The co-adjoint operator can be expressed as the block diagonal matrix,
\begin{align}
    \mathrm{ad}^*_\xi = \mathrm{diag}[0_{3\times 3}, -\hat\Omega, -\hat\Omega_R, -\hat\Omega_L, - \hat\Omega_A] \in \Re^{15\times 15}
\end{align}
since this operation on $ \SO $ is given by $ \mathrm{ad}^*_\Omega \eta = \hat{\Omega}^T \eta $.

\begin{prop}
The Euler--Lagrange equations for the flapping wing UAV are given according to \eqref{eqn:EL_G} as
\begin{gather}
	\mathbf{J}_g(\dot \xi) - \ad^*_\xi \cdot \mathbf{J}_g(\xi) + \mathbf{L}_g(\xi) \xi  = \mathbf{f}_a + \mathbf{f}_g + \mathbf{f}_\tau. \label{eqn:EL}
\end{gather}
Here, effects of the configuration dependent inertia is represented by the  matrix  $\mathbf{L}_g(\xi) = \mathbf{K}_g(\xi)  - \frac{1}{2}\mathbf{K}^T_g(\xi)\in\Re^{15\times 15}$. 
Meanwhile, $ \mathbf{f}_g = - \T^*_e\L_g \D_g U $ is the contribution of potential energy and $ \mathbf{f}_a + \mathbf{f}_\tau $ is the non-conservative external force.
\end{prop}

\section{Optimal Control}\label{sec:periodic}

In this section, we present optimal control of flapping wing UAV inspired by Monarch. 
First, the equations are reorganized such that the flapping motion can be described by wing kinematics. 
Second, we formulate an optimization to identify a periodic motion corresponding to hovering, with an additional numerical analysis to illustrate the effects of abdomen undulation. 
Next, we present an optimal control problem to stabilize the hovering flight.

\subsection{Reduced Equations}
We are interested in the global motion of the flapping wing UAV in 3-D space which is influenced by the coupled movement of wings and abdomen.
So, in this section we consider a simpler case of the dynamics wherein we prescribe the motion of wings and abdomen.
That is, we obtain equations governing the evolution of $ (x,R) $ for given functions $Q_R(t), Q_L(t), Q_A(t)$.
This is reasonable as the inertia of the wing and the abdomen are relatively small, and the corresponding torques at the joint can be reconstructed by dynamic inversion. 
This corresponds to a specific choice to formulate the maneuver with wing kinematics, and the effects of dynamic coupling are still accounted completely by \eqref{eqn:EL}.

\begin{defn}
	The configuration variables are decomposed into the free part and the prescribed part as
	\begin{gather}
		g_1 = (x, R), \quad \xi_1 = [\dot x, \Omega], \label{eqn:g1xi1}\\
		g_2 = (Q_R, Q_L, Q_A), \quad \xi_2 = [\Omega_R, \Omega_L, \Omega_A].\label{eqn:g2xi2}
	\end{gather}
	with  $g=(g_1,g_2)$ and $\xi=(\xi_1,\xi_2)$.
\end{defn}

\begin{defn}[Configuration subspaces]
	Decompose all $ 15 \times 15 $ matrices into $ \braces{(6 \times 6), (6 \times 9), (9 \times 6), (9 \times 9)} $ blocks.
	For instance, $ \mathbf{J}_g $ can be decomposed into $ \mathbf{J}_{11} \in \Re^{6 \times 6} $, $ \mathbf{J}_{21} \in \Re^{9 \times 6} $ and so on.
	Similarly a vector $ \mathbf{f} \in \Re^{15} $ can be divided into $ \mathbf{f}_1 \in \Re^{6} $ and $ \mathbf{f}_2 \in \Re^{9} $.
\end{defn}

\begin{prop}\label{prop:EL_xR}
	The derivative of $ \xi_1 $ for given $(g_2,\xi_2,\dot\xi_2)$ can be evaluated as,
	\begin{align}
		\dot \xi_1 &= (\mathbf{J}_{11}-C\mathbf{J}_{21})^{-1} \left[ (\ad^*_{\xi_1}\mathbf{J}_{11}-C\ad^*_{\xi_2} \mathbf{J}_{21} )\xi_1 \right. \nonumber\\
		& \qquad - (\mathbf{L}_{11}-C\mathbf{L}_{21})\xi_1 - (\mathbf{J}_{12}-C\mathbf{J}_{22})\dot \xi_2 \nonumber \\
		& \qquad +(\ad^*_{\xi_1}\mathbf{J}_{12}-C\ad^*_{\xi_2} \mathbf{J}_{22} )\xi_2 - (\mathbf{L}_{12}-C\mathbf{L}_{22})\xi_2 \nonumber \\ 
		&\qquad \left. + \mathbf{f}_{a_1}+\mathbf{f}_{g_1}-C(\mathbf{f}_{a_2}+\mathbf{f}_{g_2}) \right],\label{eqn:EL_xR}
	\end{align}
	where,
	\begin{align*}
		C = \begin{bmatrix} 0 & 0 & 0 \\
			-Q_R & -Q_L & -Q_A \end{bmatrix} \in \Re^{6\times 9}.
	\end{align*}
\end{prop}

\begin{pf}
The Euler--Lagrange equations~\eqref{eqn:EL} for the full configuration can be decomposed into two parts as,
\begin{align}
	\mathbf{J}_{11}\dot \xi_1 + \mathbf{J}_{12}\dot\xi_2 &-\ad^*_{\xi_1}\cdot(\mathbf{J}_{11}\xi_1 + \mathbf{J}_{12}\xi_2) + \nonumber\\
	&\mathbf{L}_{11}\xi_1 + \mathbf{L}_{12}\xi_2 = \mathbf{f}_{a_1} + \mathbf{f}_{g_1} + \mathbf{f}_{\tau_1},\label{eqn:J11dotxi_1}\\
	\mathbf{J}_{21}\dot \xi_1 + \mathbf{J}_{22}\dot\xi_2 &-\ad^*_{\xi_2}\cdot(\mathbf{J}_{21}\xi_1 + \mathbf{J}_{22}\xi_2) + \nonumber \\
	& \mathbf{L}_{21}\xi_1 + \mathbf{L}_{22}\xi_2 = \mathbf{f}_{a_2} + \mathbf{f}_{g_2} + \mathbf{f}_{\tau_2}.\label{eqn:J22dotxi_2}
\end{align}
Here the external control torques $(\tau_R,\tau_L,\tau_A)$ are unknown since we are directly specifying the wing and abdomen configuration.
So, the above two equations are coupled by these torques through the relation,
\begin{align*}
	\mathbf{f}_{\tau_1} = \begin{bmatrix} 0 & 0 & 0 \\
		-Q_R & -Q_L & -Q_A \end{bmatrix} \mathbf{f}_{\tau_2} 
	= C \mathbf{f}_{\tau_2},
\end{align*}
from \eqref{eqn:f_tau}.
To remove these terms, we calculate \eqref{eqn:J11dotxi_1} - $ C \times $ \eqref{eqn:J22dotxi_2} to obtain
\begin{align*}
	&(\mathbf{J}_{11}-C\mathbf{J}_{21})\dot \xi_1 -(\ad^*_{\xi_1}\mathbf{J}_{11}-C\ad^*_{\xi_2} \mathbf{J}_{21} )\xi_1 + (\mathbf{L}_{11}-C\mathbf{L}_{21})\xi_1 \nonumber \\ 
	& \qquad =  - (\mathbf{J}_{12}-C\mathbf{J}_{22})\dot \xi_2 +(\ad^*_{\xi_1}\mathbf{J}_{12}-C\ad^*_{\xi_2} \mathbf{J}_{22} )\xi_2  \nonumber \\ 
	& \qquad \qquad - (\mathbf{L}_{12}-C\mathbf{L}_{22})\xi_2 + \mathbf{f}_{a_1}+\mathbf{f}_{g_1}-C(\mathbf{f}_{a_2}+\mathbf{f}_{g_2})
\end{align*}
which is rearranged into the equation in \eqref{eqn:EL_xR}.

The control toques $(\tau_R,\tau_L,\tau_A)$ necessary to specify motion of the wings and abdomen can then be obtained from \eqref{eqn:J22dotxi_2} by substituting the integrated $(g_1,\xi_1)$.
\hspace*{\fill} \qed
\end{pf}

\subsection{Wing and Abdomen Kinematics}

Since we are prescribing the second set of configuration in \eqref{eqn:g2xi2}, it would be simpler to parameterize the trajectories of these variables.
Consider the model utilized in~\cite{tejaswi2020effects} for the motion of the wing relative to the body.
Let $f\in\Re$ be the flapping frequency in $\mathrm{Hz}$ and $T=\frac{1}{f}$ be the corresponding time period in seconds. 

\begin{itemize}
	\item The flapping angle is parameterized as,
	\begin{align}
		\phi(t) & = \frac{\phi_m}{\sin^{-1} \phi_K}\sin^{-1}(\phi_K\cos(2\pi f t)) + \phi_0,\label{eqn:phi}
	\end{align}
	where $\phi_m\in\Re$ is the amplitude, $\phi_0\in\Re$ is the offset, and $0 < \phi_K \leq 1$ determines waveform shape.
	\item The pitch angle is given by,
	\begin{align}
		\theta(t) = \frac{\theta_m}{\tanh \theta_C} \tanh( \theta_C \sin(2\pi f t + \theta_a)) +\theta_0,\label{eqn:theta}
	\end{align}
	where $\theta_m\in\Re$ is the amplitude of pitching, $\theta_0\in\Re$ is the offset, $\theta_C\in(0,\infty)$ determines the waveform, and $\theta_a\in(-\pi,\pi)$ describes phase offset. 
	\item Finally, the deviation angle is given by 
	\begin{align}
		\psi(t) = \psi_m \cos(2\pi \psi_N f t + \psi_a) + \psi_0,\label{eqn:psi}
	\end{align}
	where $\psi_m\in\Re$ is the amplitude, $\psi_0\in\Re$ is the offset, and the parameter $\psi_a\in(-\pi,\pi)$ is the phase offset. 
\end{itemize}
Using these Euler angles, the attitude, angular velocity and acceleration of the wings can be constructed.

Next, the attitude of the abdomen relative to the body can be considered as $Q_A(t) = \exp(\theta_A(t)\hat e_2)$.
This is motivated by the flight characteristics of a live Monarch butterfly which exhibits a nontrivial pitching motion of the abdomen (see ~\cite{sridhar2020geometric}).
Here, the relative pitch angle is taken to be
$
\theta_A(t) = \theta_{A_m} \cos{(2 \pi f t + \theta_{A_a})} + \theta_{A_0},
$
for fixed parameters $\theta_{A_m},\theta_{A_a},\theta_{A_0}\in\Re$.

\subsection{Periodic Motion}

\begin{figure}
	\centerline{
		\subfigure[Position of body $x$]{
			\includegraphics[width=0.5\linewidth]{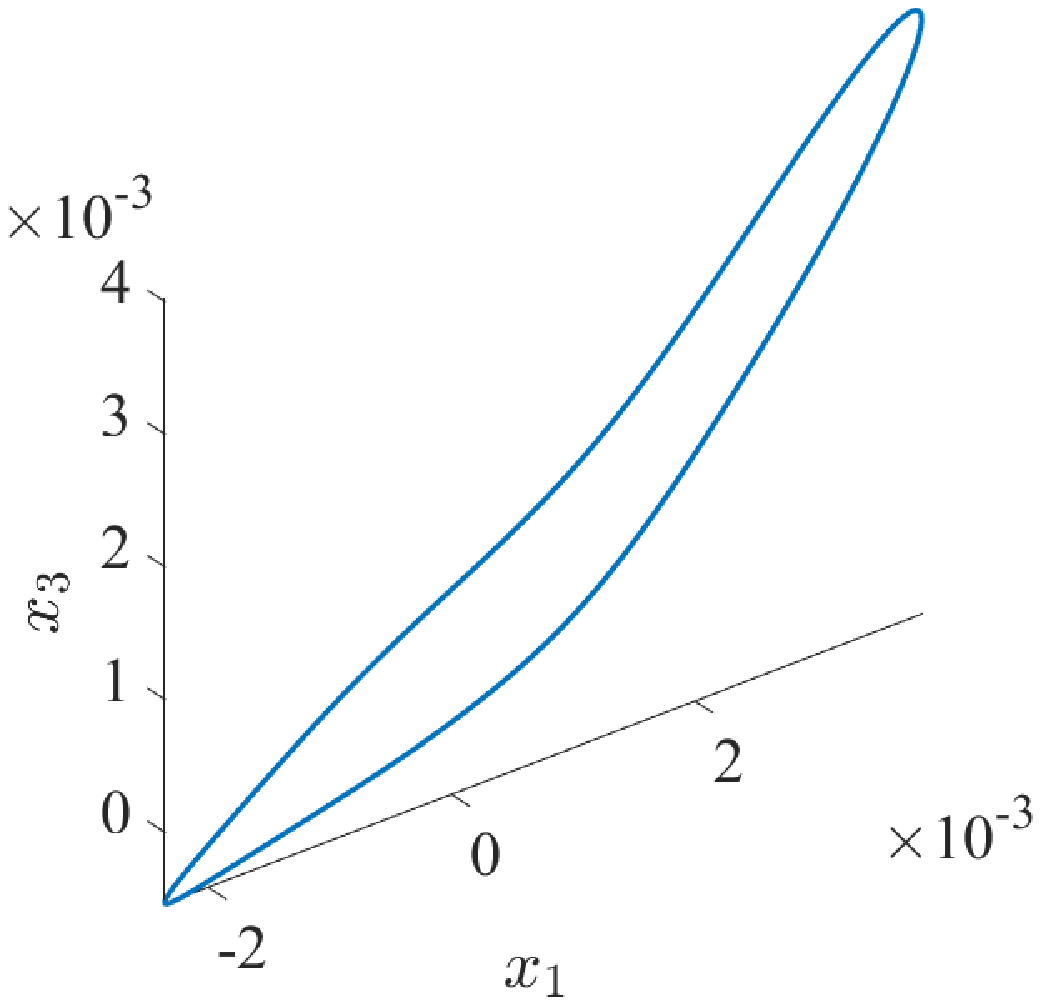}
		}
		\hfill
		\subfigure[Velocity of body $\dot x$]{
			\includegraphics[width=0.5\linewidth]{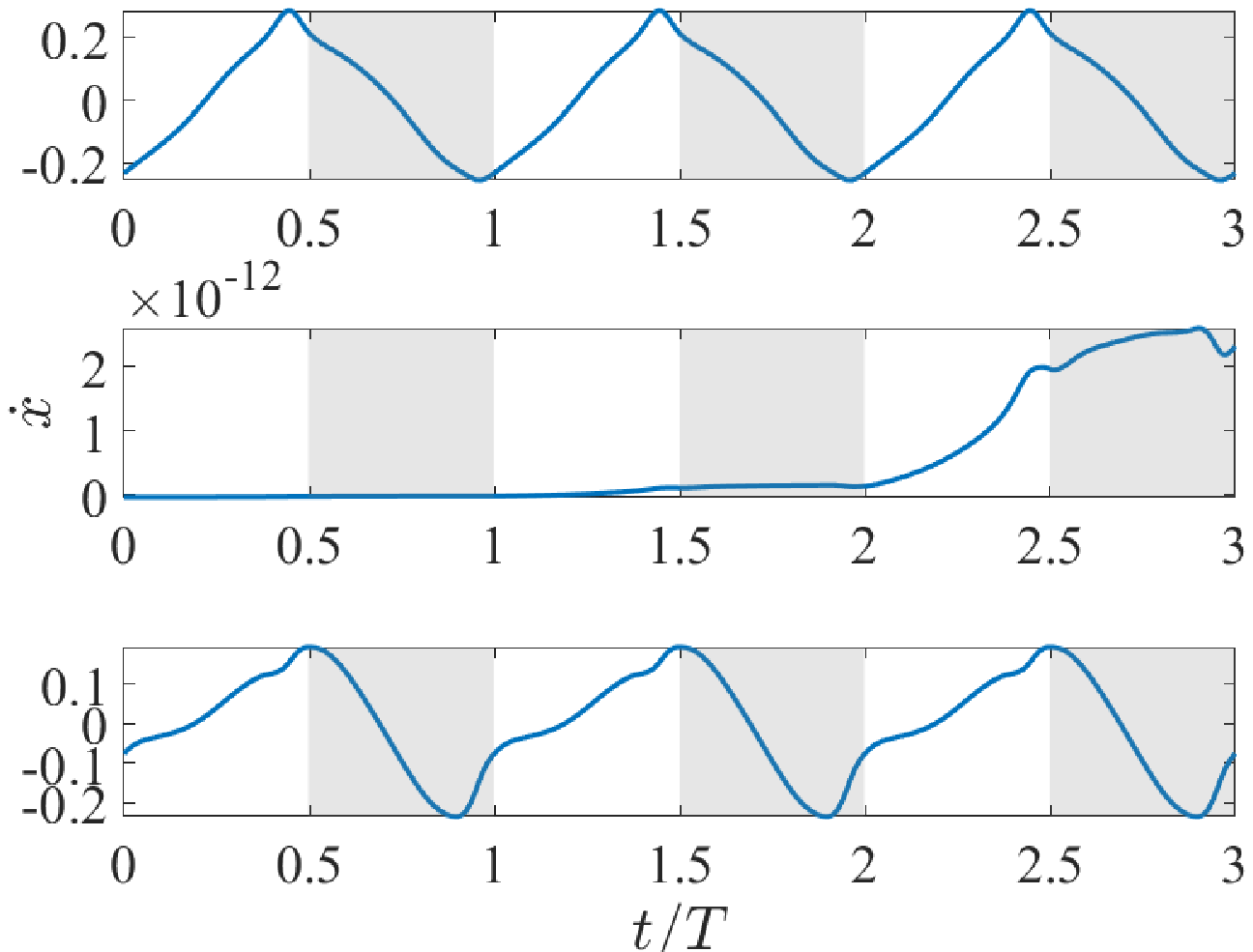}
		}
	}
	\centerline{
		\subfigure[Body pitch (in degrees) and angular velocity along 2nd axis]{
			\includegraphics[width=0.5\linewidth]{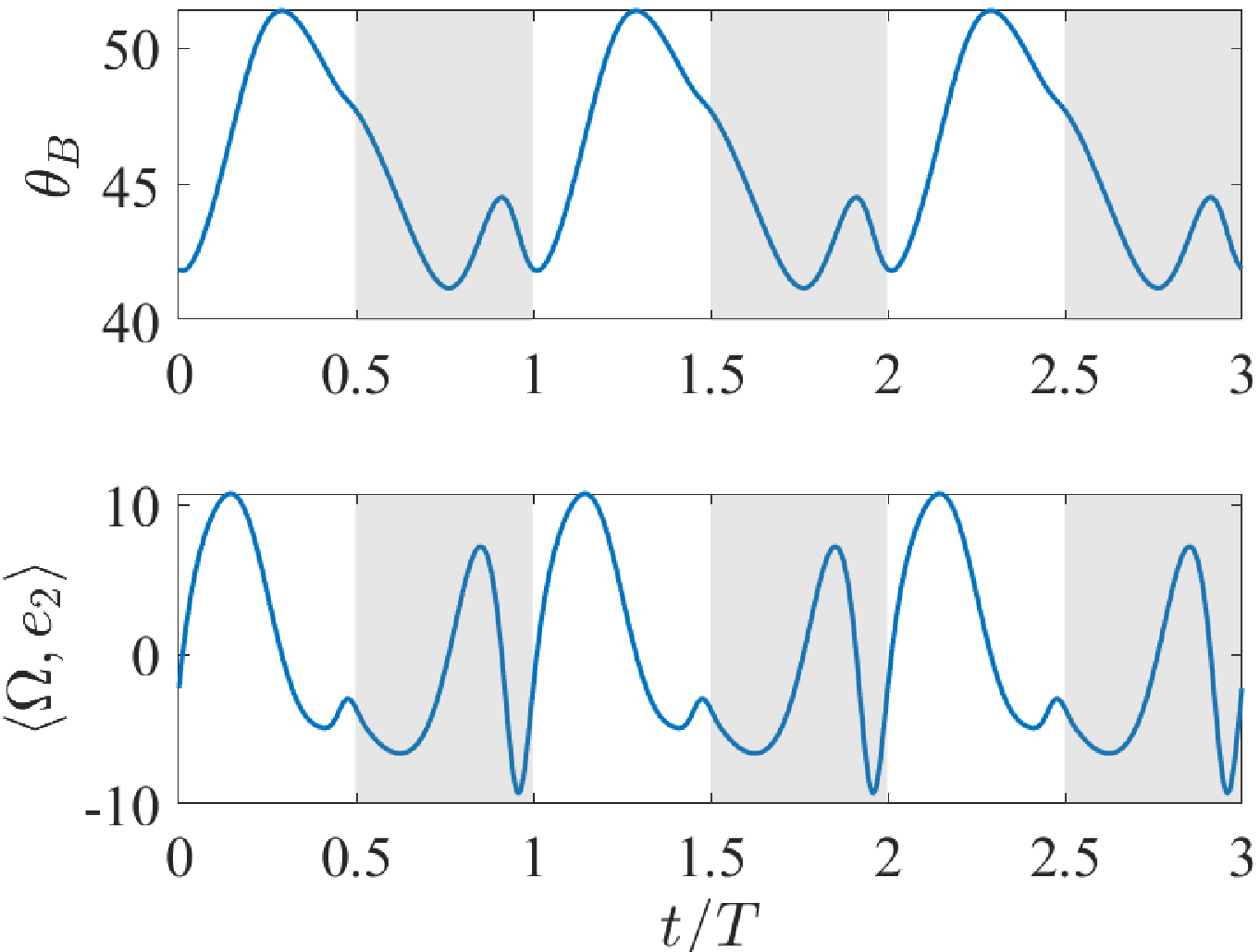}
		}
		\hfill
		\subfigure[Prescribed wing kinematics and abdomen undulation (in degrees)]{
			\includegraphics[width=0.5\linewidth]{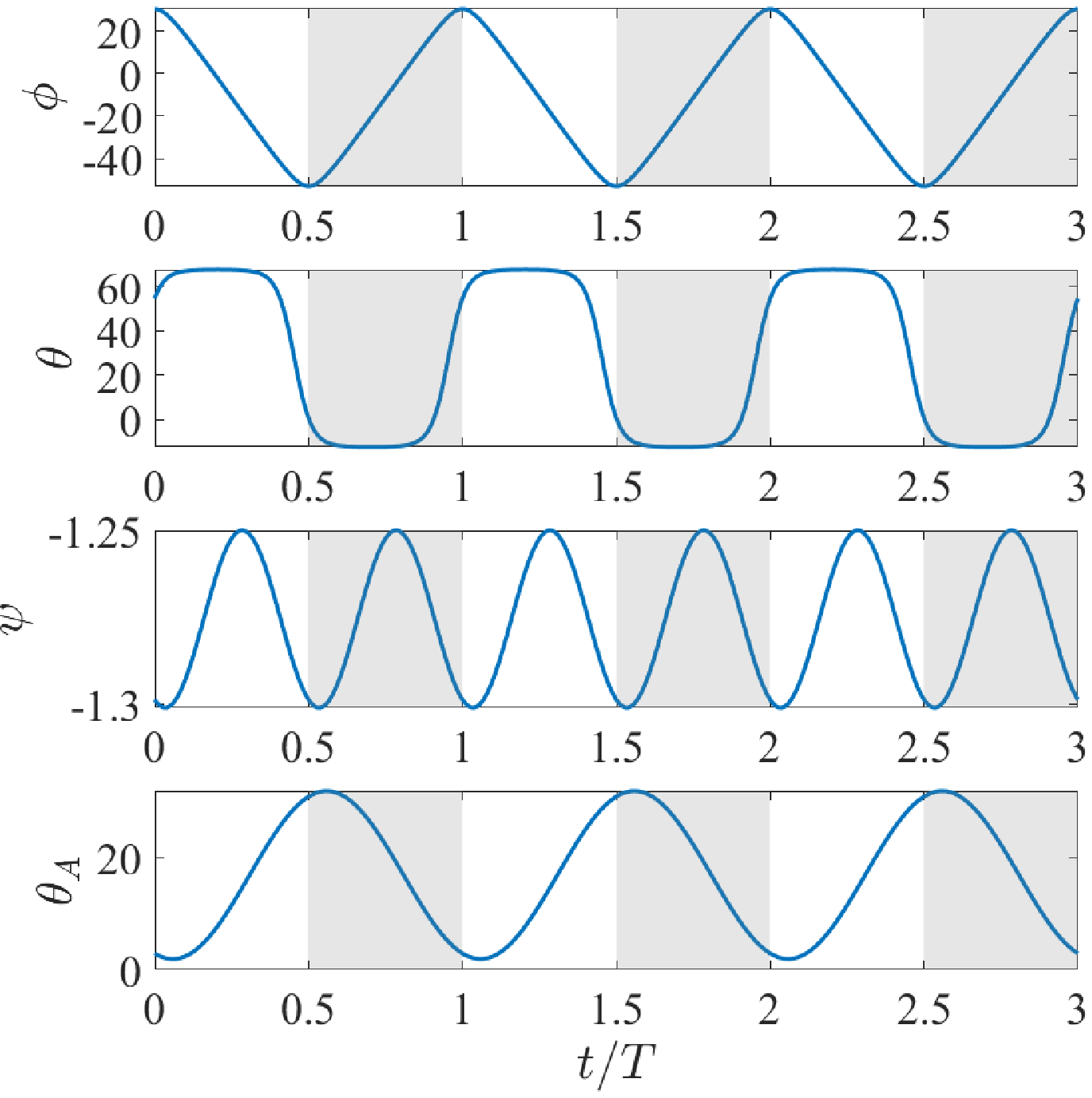}
		}
	}
	\caption{Hovering periodic orbit generated using optimized parameters; shaded region corresponds to downstrokes}
	\label{fig:hover_pos_vel}
\end{figure}

The above dynamic model yields the position and the attitude trajectory of the body for given kinematics of wings and abdomen.
We first need to construct the kinematics of wings and abdomen for a particular maneuver. 
This is challenging due to the complexities of the dynamics and the relatively large number of free parameters in the wing kinematics.
Here we focus on the case of hovering flight, where the position and the attitude returns to the initial value after each flapping period. 
This result can be easily extended to other maneuvers such as forward flight or climbing.

This is addressed by a constrained optimization to minimize a performance index while ensuring that the motion is periodic.
The parameters being optimized over characterize FWUAV wing kinematics and abdomen undulation along with the initial conditions.
More specifically, this is formulated as follows.
\begin{itemize}
    \item The objective function is 
        \begin{equation}\label{eqn:obj_func}
            J = w_1 \int_{0}^{T} |E(t)| dt + w_2 \int_{0}^{T} |\dot{E}(t)| dt,
        \end{equation}
        where $w_1,w_2\in\Re$, and
        $
            E(t) = \frac{1}{2} m \norm{\dot x(t)}^2 - mge_3^T x(t)
        $
        is the sum of the kinetic energy and the gravitational potential energy.
        This is to minimize the variation of the energy while penalizing abrupt changes. 
    \item The optimization parameters are given by
        \begin{itemize}
            \item flapping frequency: $f$ and stroke plane angle: $\beta$
            \item wing kinematics: $(\phi_m, \phi_K, \phi_0)$, $(\theta_m, \theta_C, \theta_0, \theta_a)$, $(\psi_m, \psi_0, \psi_a)$
            \item abdomen undulation : $(\theta_{A_m}, \theta_{A_0}, \theta_{A_a})$
            \item initial translational velocity: $\dot x(0) \in \Re^3$
            \item initial attitude, angular velocity along 2nd axis: $ \theta_B(0)\ s.t.\ R(0) = \exp(\theta_B(0) \hat{e}_2) $, $ \Omega_2(0) = \pair{\Omega, e_2} $
        \end{itemize}
    \item We ensure periodic motion by imposing the constraints:
        $
            x(0) = x(T),\quad  \dot{x}(0) = \dot{x}(T).
        $
        Furthermore, there are additional constraints to avoid physically infeasible flapping,
        $
            |\phi_m| + |\phi_0| < \pi/2,
        $
        along with prescribed hard bounds on other parameters.
        It is also assumed that the motion of wings is symmetric to each other in this simple maneuver.
        \item The physical properties of the FWUAV including the wing morphological parameters like $ J_i, \mu_i $ are taken to be similar to those of an actual Monarch.
        Their specific values are given in~\cite{sridhar2020geometric}.
        \item The aerodynamic properties including lift and drag coefficients are adopted from experimental data in~\cite{dickinson1999wing,sane2001control}.
        \cite{tejaswi2020effects} presents these expressions along with their relations to the actual aerodynamic forces and torques in \eqref{eqn:f_a}.
        Furthermore, only the wings are assumed to generate aerodynamic forces since the projected area of the body and abdomen is not significant.
\end{itemize}

This problem is solved via global optimization techniques such as \texttt{multistart} in MATLAB.
The corresponding optimized parameters are summarized in Table \ref{tab:hover_params}, and the resulting maneuver is illustrated in Figure \ref{fig:hover_pos_vel}.
Note that since this maneuver is in the $x$-$z$ plane, the relative attitude and angular velocity of the body are non-zero only along the $ y $ axis as shown in Figure \ref{fig:hover_pos_vel}.(c).
Compared with~\cite{tejaswi2020effects} where the periodic orbit is constructed for the translational dynamics, this provides the periodic motion for the coupled translational and rotational motion in the higher-dimensional space.

\begin{table}[h!]
    \newcolumntype{m}{>$l<$} 
	\caption{Optimized parameters}\label{tab:hover_params}
	\begin{center}
		\begin{tabular}{|m|m|m|}
            \hline
            \text{Parameters} & \text{With abdomen} & \text{Without abdomen}\\
            & \text{undulation} & \text{undulation} \\\hline
            f & 11.7575 & 11.3975 \\
            \beta & -0.0087 & 0.2014 \\
            \phi_m & 0.7271 & 0.6655 \\
            \phi_K & 0.9493 & 0.0138 \\
            \phi_0 & -0.1977 & -0.0434 \\
            \theta_m & 0.6981 & 0.6980 \\
            \theta_C & 2.8289 & 2.9968 \\
            \theta_0 & 0.4843 & 0.3503 \\
            \theta_a & 0.2905 & 0.3971 \\
            \psi_m & 0.0004 & 0.0003 \\
            \psi_N & 2 & 2 \\
            \psi_0 & -0.0223 & -0.0400 \\
            \psi_a & 2.7130 & 3.1109 \\
            \theta_{A_m} & 0.2618 & \text{\textemdash\textemdash} \\
            \theta_{A_0} & 0.2950 & 0.7667 \\
            \theta_{A_a} & 2.7743 & \text{\textemdash\textemdash} \\
            \dot x_1(0) & -0.2332 & -0.2437 \\
            \dot x_2(0) & 0.0000 & 0.0000 \\
            \dot x_3(0) & -0.0764 & -0.0859 \\
            \theta_B(0) & 0.7314 & 0.5666 \\
            \Omega_2(0) & -2.2583 & -0.1709 \\
            \hline
            \text{Optimized}\ J & 0.0787 & 0.0890 \\\hline
        \end{tabular}\\[0.1cm]
        ($ f_{natural} = 10.2247\,\si{Hz}, \psi_N = 2 $)
    \end{center}
\end{table}

\subsection{Effects of Abdomen}

Now we study the influence of abdomen undulation on the periodic maneuver and the performance index.
As a comparison, we identify another periodic orbit assuming that the abdomen is at a fixed relative attitude with respect to the body.
The second column of Table \ref{tab:hover_params} lists the optimized parameters wherein $ \theta_{A_0} $ is the constant relative pitch of abdomen.

\begin{figure}
	\centerline{
		\subfigure[Energy]{
			\includegraphics[width=0.5\linewidth]{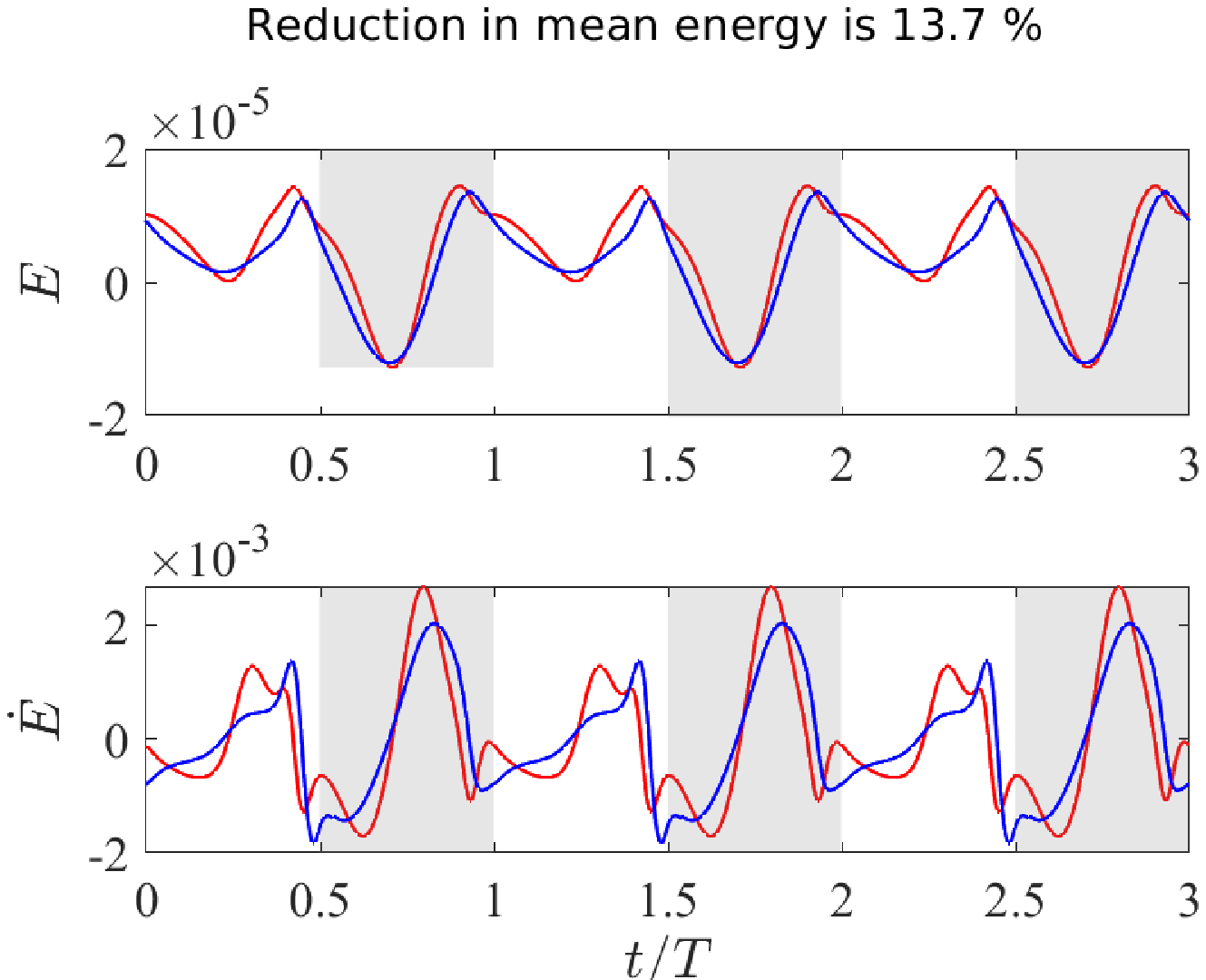}
		}
		\hfill
		\subfigure[Power]{
			\includegraphics[width=0.5\linewidth]{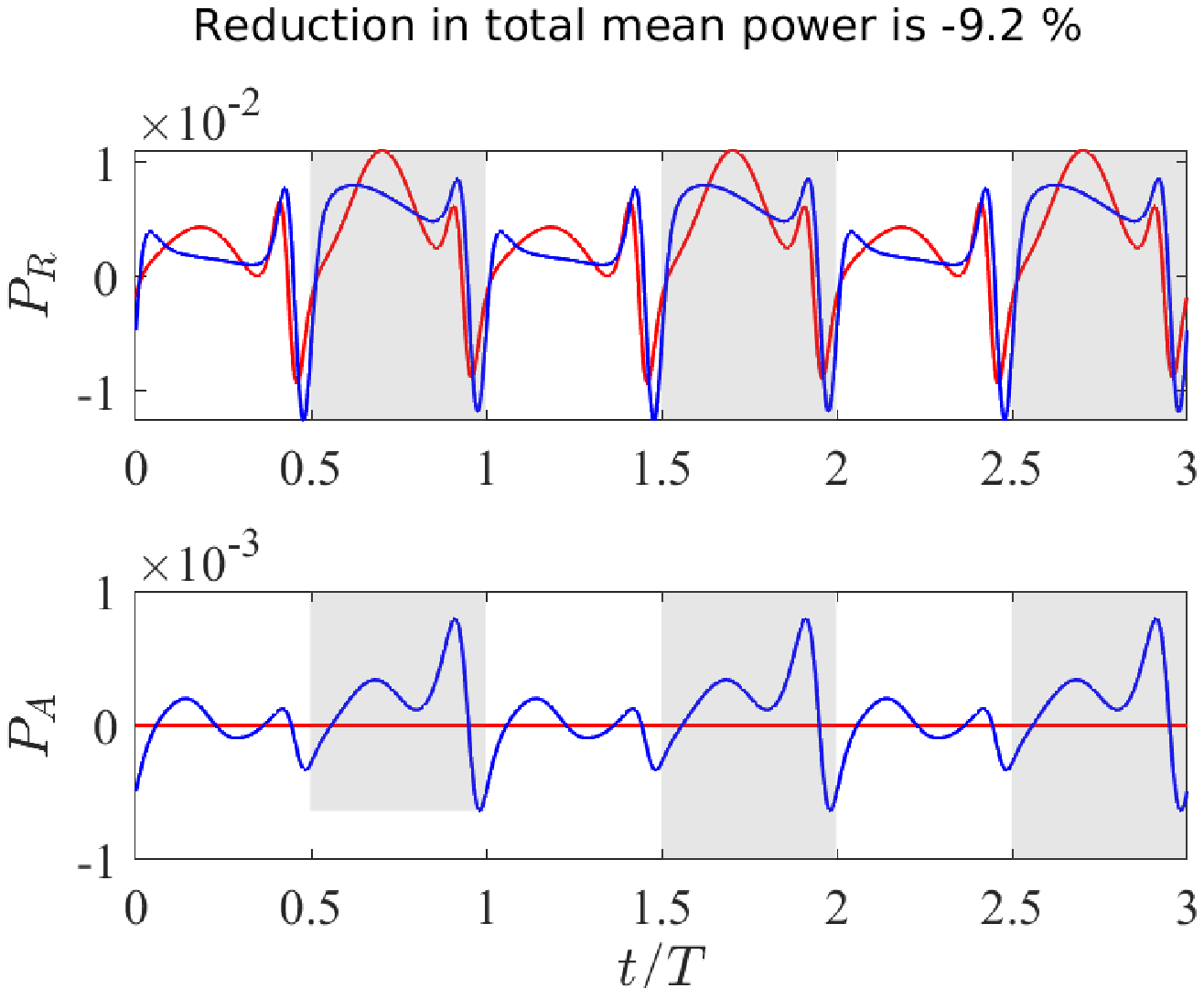}
		}
	}
	\centerline{
		\subfigure[Torque]{
			\includegraphics[width=0.5\linewidth]{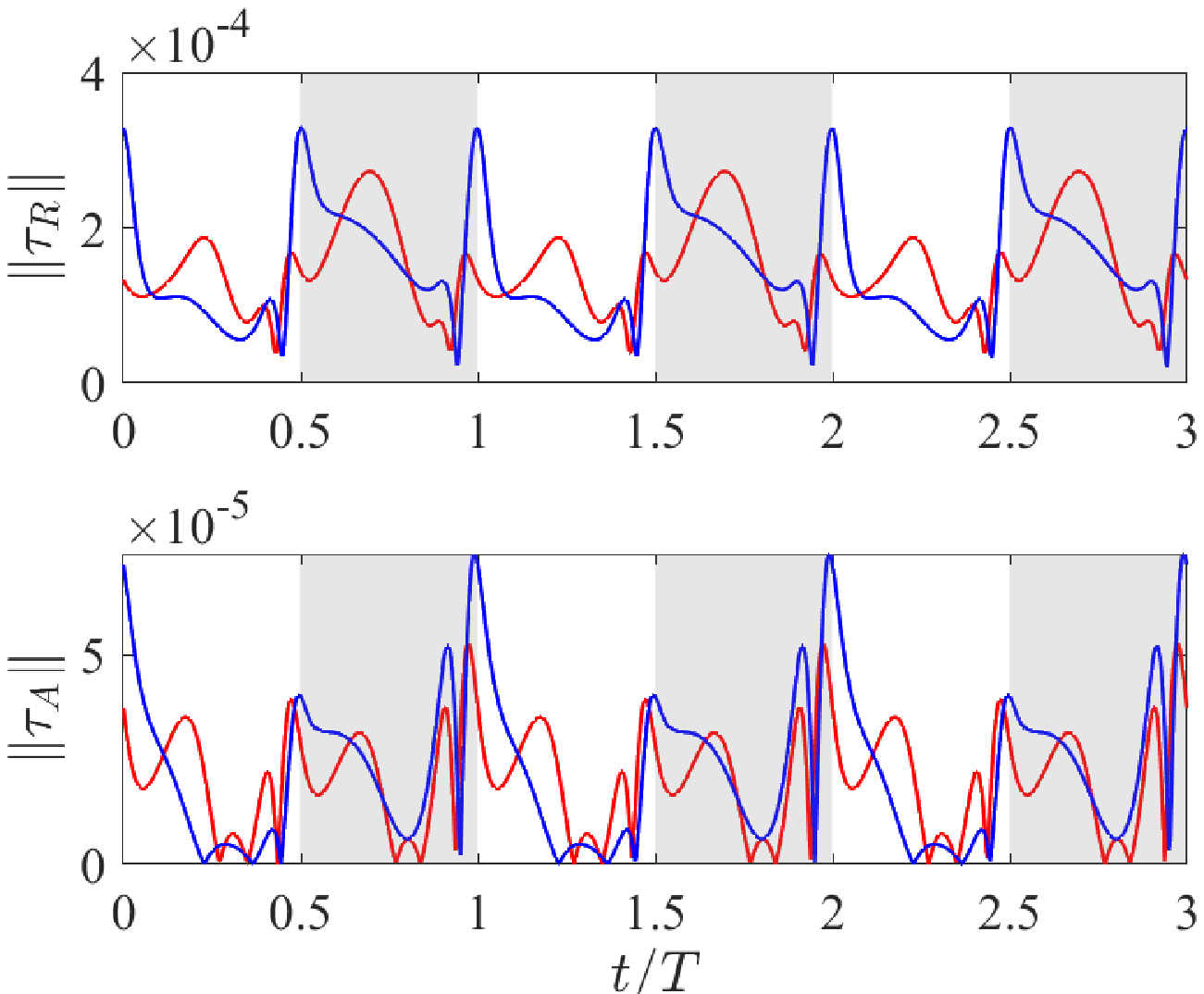}
		}
	}
    \caption{Comparison between hovering with abdomen undulation (blue) and hovering without abdomen undulation (red)
	}
	\label{fig:hover_comp_ab}
\end{figure}

We observe that the objective function is decreased by about $12\%$ when there is abdomen undulation when compared to no abdomen undulation.
Since $ J $ in \eqref{eqn:obj_func} is composed of energy and its derivative, they are also reduced in the case of abdomen undulation as seen in Figure \ref{fig:hover_comp_ab}.(a).
This is not surprising as there are additional degrees of freedom that are utilized to minimize the objective function further.

Finally, the dynamical equations are utilized to obtain the control torques $(\tau_R,\tau_L,\tau_A)$ as shown in the proof of Proposition \ref{prop:EL_xR}.
For this numerical experiment, their magnitudes are illustrated at Figure \ref{fig:hover_comp_ab}.(c).
Also, the power due to these external torques can be calculated as $P_R = \tau_R^T (Q_R \Omega_R)$ and $P_A = \tau_A^T (Q_A \Omega_A)$ in the body frame.
Figure \ref{fig:hover_comp_ab}.(b) compares these values for the cases with and without abdomen undulation.
Note that since the wings move symmetrically, $ \norm{\tau_R} = \norm{\tau_L} $ and $ P_R = P_L $.


\subsection{Optimal Control}

Now we formulate an optimal control problem such that an arbitrary trajectory asymptotically converges to the optimized periodic orbit for the hovering flight that we have obtained in the previous section. 
More specifically, let $\mathbf{x}(t) = (g_1(t),\xi_1(t)) =  (x(t), R(t), \dot x(t), \Omega(t)) $ represent the state of the FWUAV for the translational and rotational dynamics of the body.
We have already obtained a periodic reference trajectory $ \mathbf{x}_d(t) = (x_d(t), R_d(t), \dot x_d(t), \Omega_d(t)) $.
The objective is to adjust the control parameters such that $ \mathbf{x}(t) \to \mathbf{x}_d(t)$.

We have various parameters in the definition of the wing kinematics in \eqref{eqn:phi}--\eqref{eqn:psi}.
Instead of numerically optimizing all of those parameters by brute-force, we identity a smaller set of parameters by investigating the effects of those on aerodynamic forces.
So we choose $ N_\Delta = 6 $ specific control parameters:
\begin{equation}\label{eqn:params}
	\Delta = [\Delta \phi_{m_s}, \Delta\theta_{0_s}, \Delta \phi_{m_k}, \Delta\phi_{0_s}, \Delta\theta_{0_k}, \Delta\psi_{0_k}],
\end{equation}
They are composed of two types:
\begin{itemize}
	\item Symmetric parameters: for instance, $ \Delta \phi_{m_s} = (\Delta \phi_{m,R} + \Delta \phi_{m,L})/2 $ which is the average change of amplitude of the flapping angle of both wings
	\item Anti-symmetric parameters: e.g., $ \Delta\phi_{m_k} = (\Delta \phi_{m,R} - \Delta \phi_{m,L})/2 $ which is the difference of flapping amplitude changes leading to a lateral force
\end{itemize}
Here,  $ \Delta \phi_{m,R} = \phi_{m,R}(t) - \phi_{m,R,d} $, i.e., the change of flapping amplitude of the right wing from the desired trajectory.
Other variables are defined similarly. 
The effects of these control parameters on the resultant force and moment are summarized as follows.
\begin{table}[h!]
	\newcolumntype{m}{>$c<$} 
	\caption{Change in average forces/moments studied near the ideal hover trajectory}
	\begin{center}
		\begin{tabular}{|m|mmmmmm|}
			\hline
			& \Delta \phi_{m_s} & \Delta\theta_{0_s} & \Delta \phi_{m_k} & \Delta\phi_{0_s} & \Delta\theta_{0_k} & \Delta\psi_{0_k} \\\hline
			\Delta \bar f_{a_1} \times 10^4 & -23  & -31 & 0 & -37 & 0 & 0 \\ 
			\Delta \bar f_{a_2} \times 10^4 & 0 & 0 & 87 & 0 & -77 & 0\\
			\Delta \bar f_{a_3} \times 10^4 & -78 & 42 & 0 & 41 & 0 & 0\\\hline
			\Delta \bar M_{a_1} \times 10^5 & 0 & 0 & 79 & 0 & 0 & 0\\
			\Delta \bar M_{a_2} \times 10^5 & -7 & 4 & 0 & 6 & 0 & 0 \\
			\Delta \bar M_{a_3} \times 10^5 & 0 & 0 & -111 & 0 & 35 & -14\\\hline
		\end{tabular}
	\end{center}
\end{table}
The proposed control parameters improve the efficiency of optimization, and the corresponding optimized trajectories are more suitable to be generalized into other maneuvers. 

To represent the variation of these control parameters over time, the flapping period $[0, T] $ divided into $N_s = 10 $ steps at which the values of the control parameters are specified.
Considering that the desired trajectory is periodic, we impose an additional constraint $ \Delta(0) = \Delta(T) = 0 $.
The value of $\Delta(t)$ between discrete steps are obtained by a piecewise linear interpolation.

The objective function is the weighted sum of the discrepancy between the desired trajectory and the controlled trajectory given by
\begin{align}
    J = \sum_{i=1}^{N_p} W_i \sqrt{\sum_j (W_{\mathbf{x}_j} (\mathbf{x}_j(t_i) - \mathbf{x}_{d_j} (t_i)))^2}
\end{align} 
where $ t_i = i \times T/N_s $.
The inner sum represents errors in the states $ \mathbf{x} $ at time $ t_i $ weighed by a factor $ W_{\mathbf{x}_j}$.
The outer sum combines the state errors over each prediction horizon $N_p$ weighed by another factor $ W_i $.
The weighting factor for the state, $ W_\mathbf{x} $ is designed to ensure that each component is scaled by its own physical characteristics.
And the weighting factor for time, $ W_i $  gradually increases over $ i $ so that the terminal state error has more weight.

We follow the formulation of model predictive control, where the prediction horizon corresponds to two flapping period, i.e., $N_p = 2N_s = 20$.
The optimization is repeated at every period to find the optimal control parameters over the prediction horizon, resulting in $120$ optimal control parameters for two periods.
Among those, the control parameters corresponding to the first period is actually implemented, and at the end of the period, optimization is repeated. 


\begin{figure}
	\centerline{
		\subfigure[Linear velocity error]{
			\includegraphics[width=0.5\linewidth]{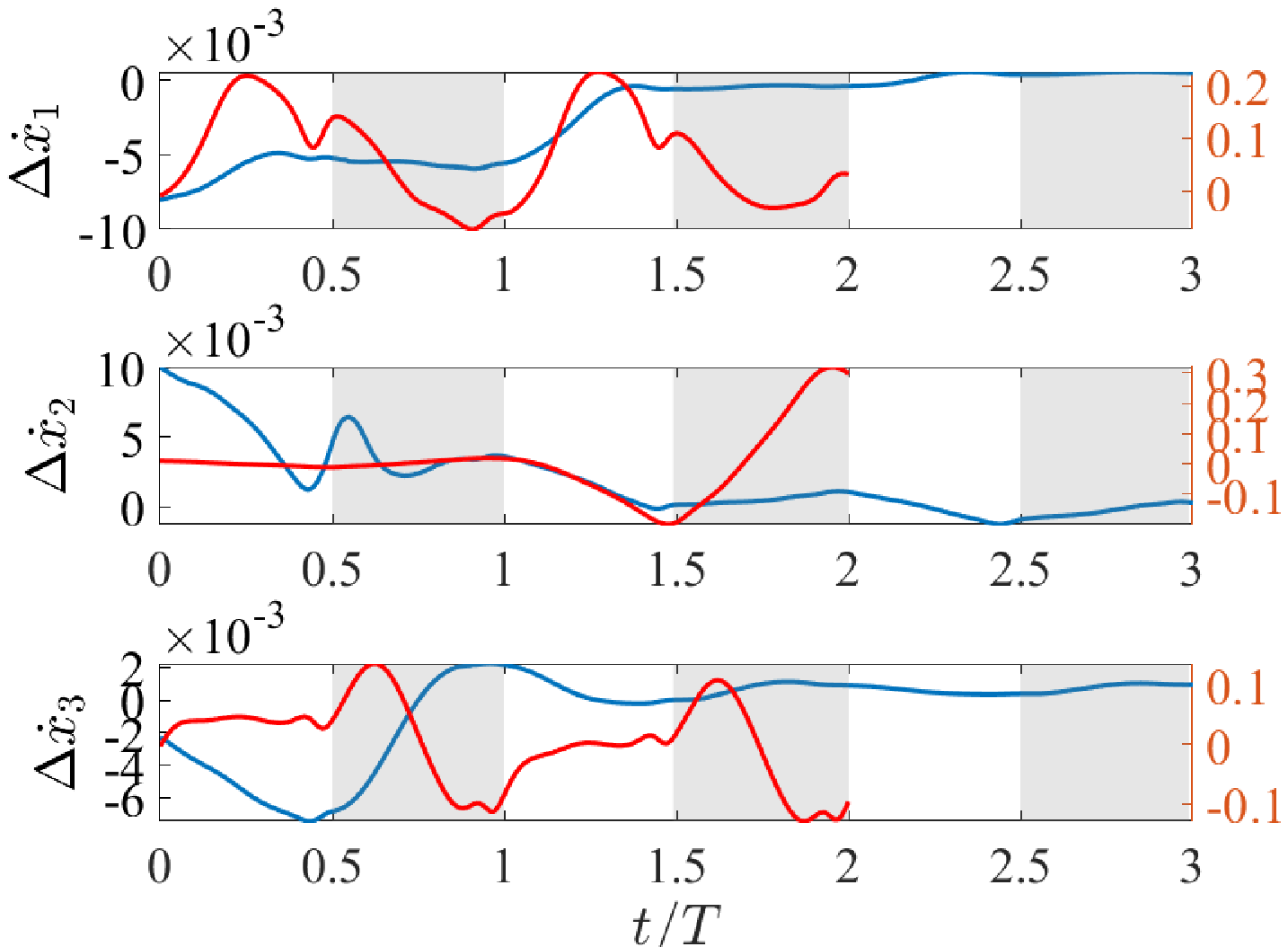}
		}
		\hfill
		\subfigure[Angular velocity error]{
			\includegraphics[width=0.5\linewidth]{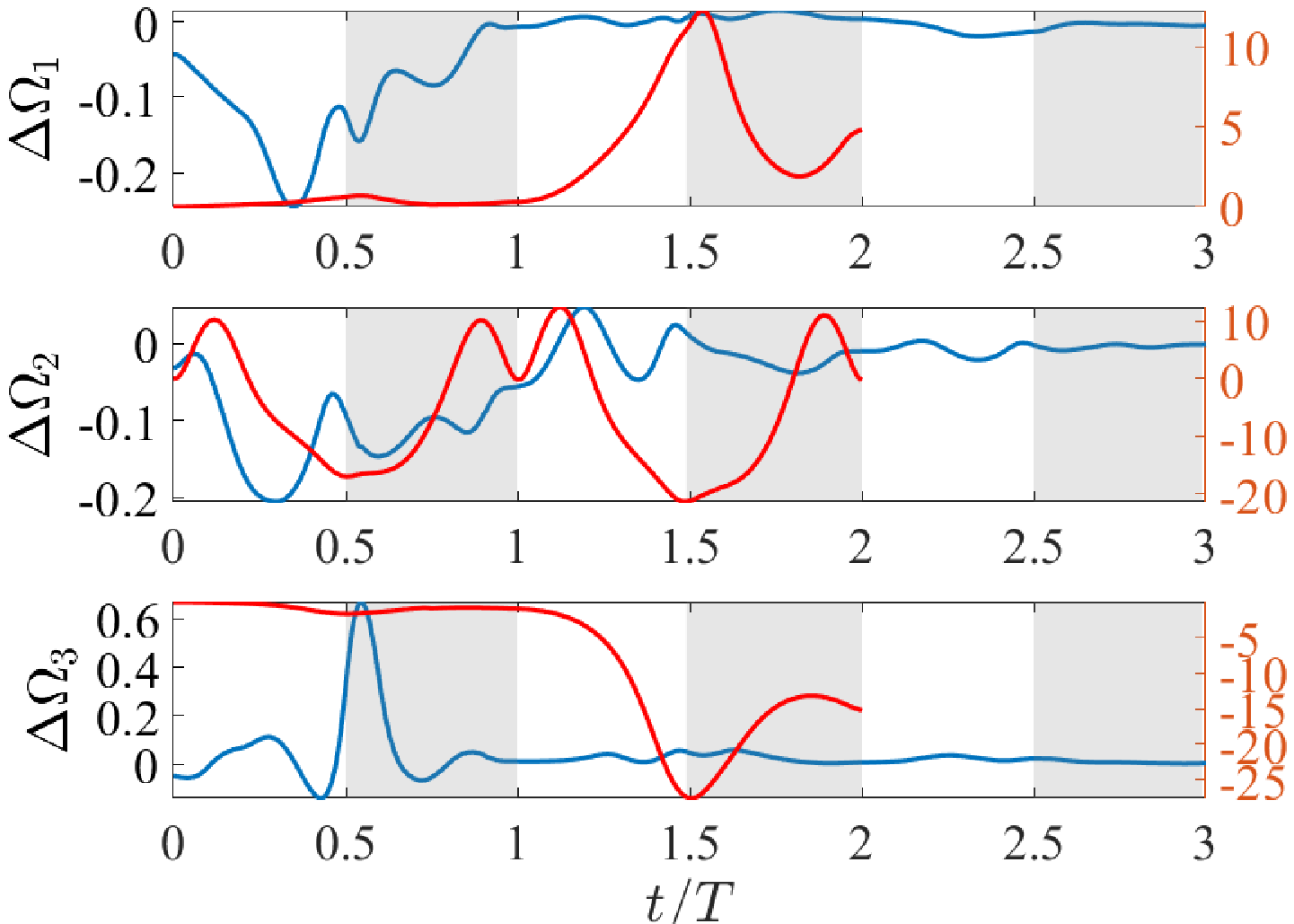}
		}
	}
	\centerline{
		\subfigure[Position and attitude error]{
			\includegraphics[width=0.5\linewidth]{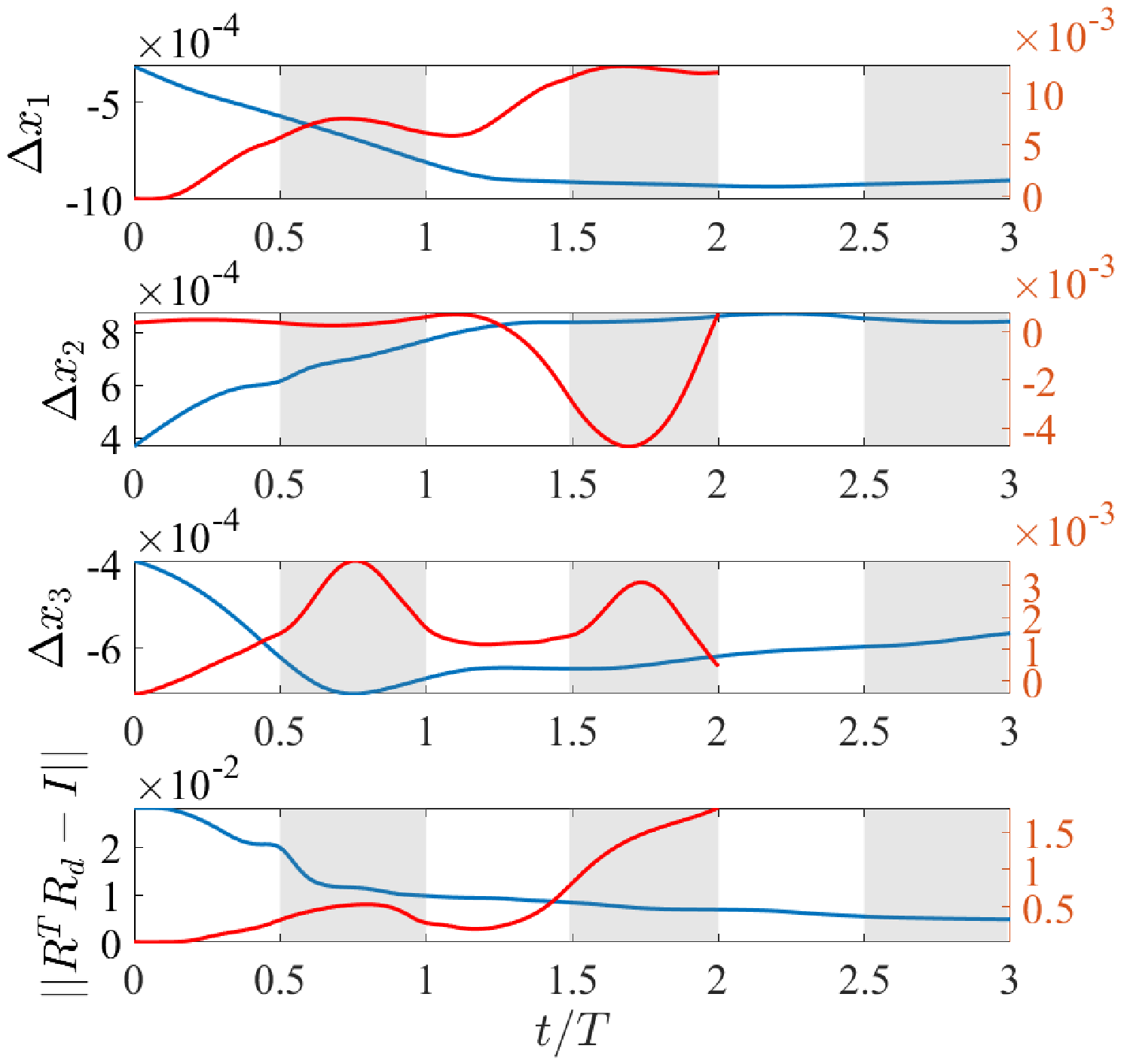}
		}
		\subfigure[Control inputs]{
			\includegraphics[width=0.475\linewidth]{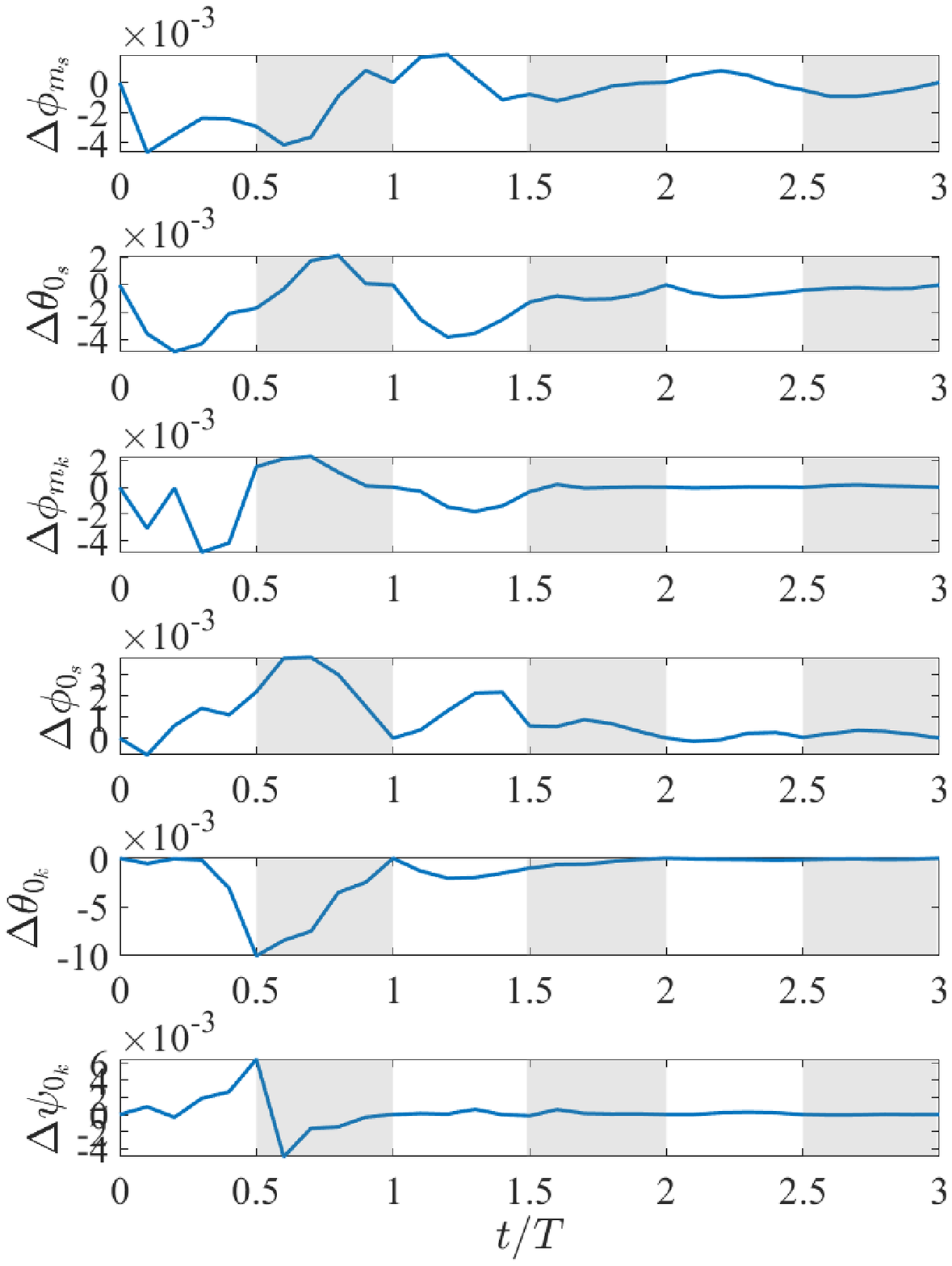}
		}
	}
	\caption{
	Optimal trajectory errors in blue with scale on the left y-axis; uncontrolled trajectory errors in red with labels on the right y-axis}
	\label{fig:hover_control}
\end{figure}

\begin{figure}
	\centerline{
		\subfigure[$ t=T $]{
			\includegraphics[trim={11cm 4cm 12cm 4cm},clip,width=0.25\linewidth]{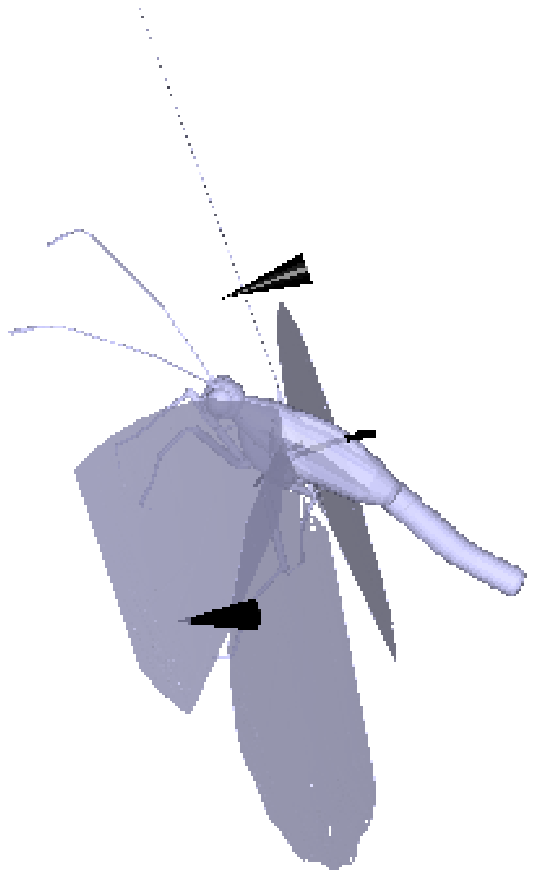}
		}
		\hfill
		\subfigure[$ t=4T/3 $]{
			\includegraphics[trim={12cm 5cm 11cm 4cm},clip,width=0.25\linewidth]{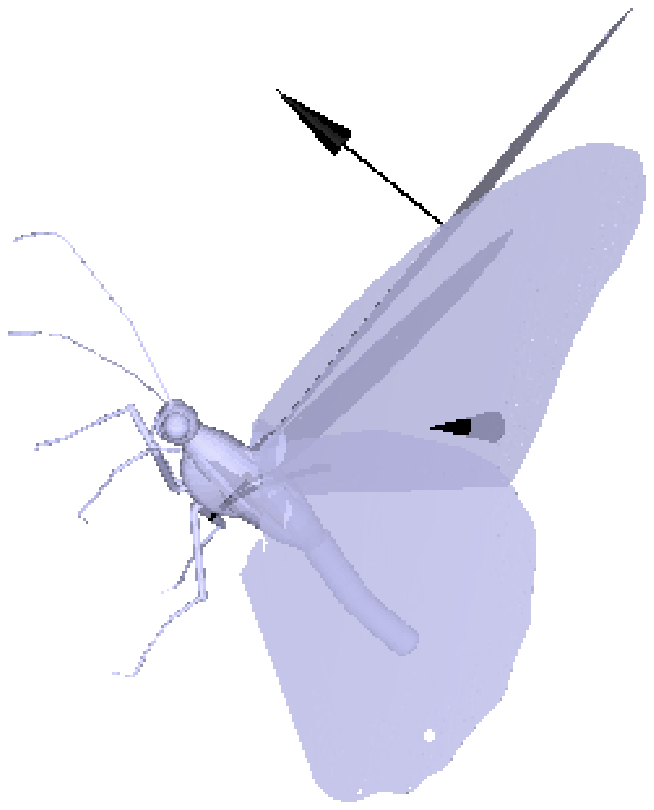}
		}
		\hfill
		\subfigure[$ t=5T/3 $]{
			\includegraphics[trim={12cm 6cm 11cm 3cm},clip,width=0.25\linewidth]{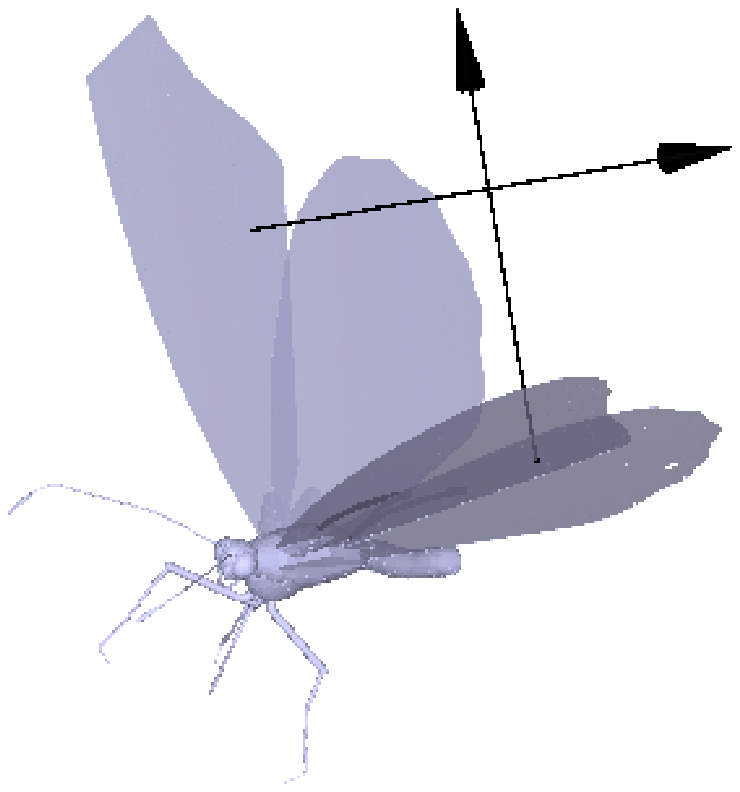}
		}
		\hfill
		\subfigure[$ t=2T $]{
			\includegraphics[trim={11cm 5cm 11cm 4cm},clip,width=0.25\linewidth]{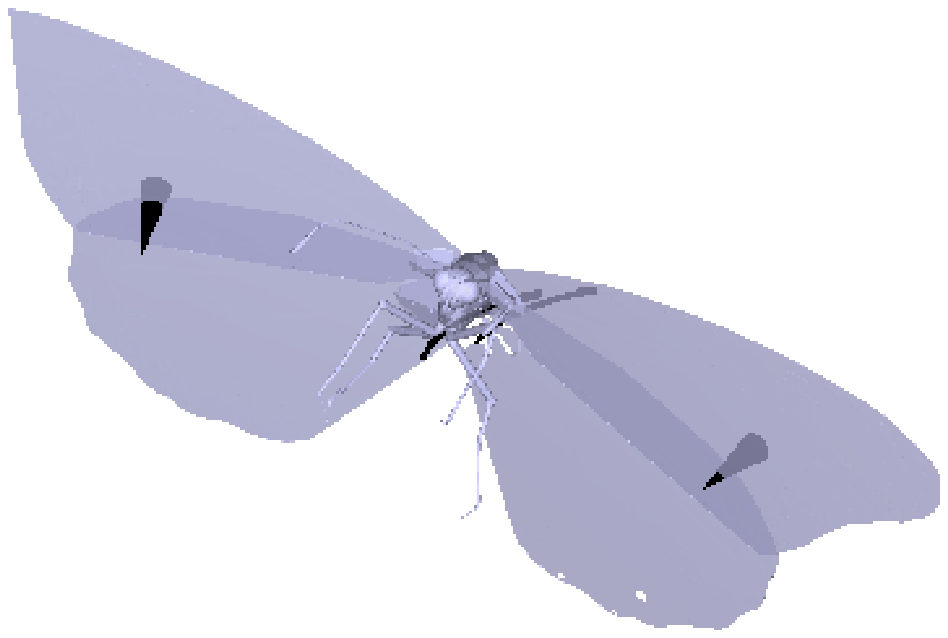}
		}
	}
	\centerline{
		\subfigure[$ t=T $]{
			\includegraphics[trim={11cm 4cm 12cm 4cm},clip,width=0.25\linewidth]{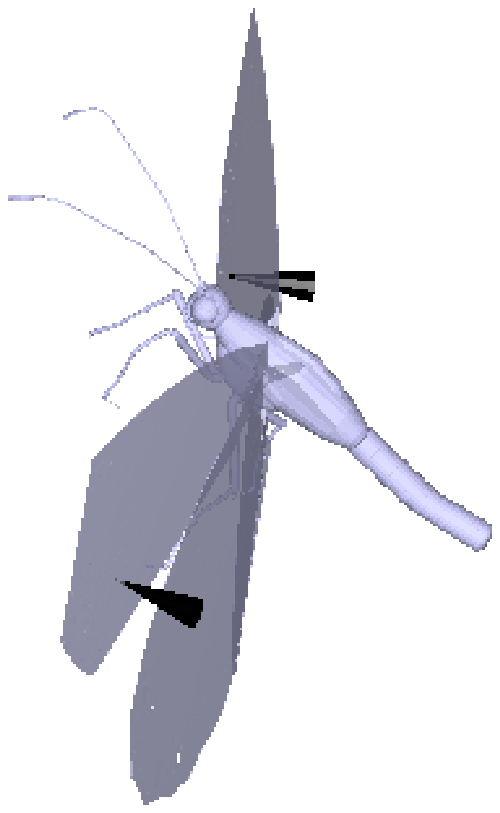}
		}
		\hfill
		\subfigure[$ t=4T/3 $]{
			\includegraphics[trim={12cm 5cm 12cm 4cm},clip,width=0.25\linewidth]{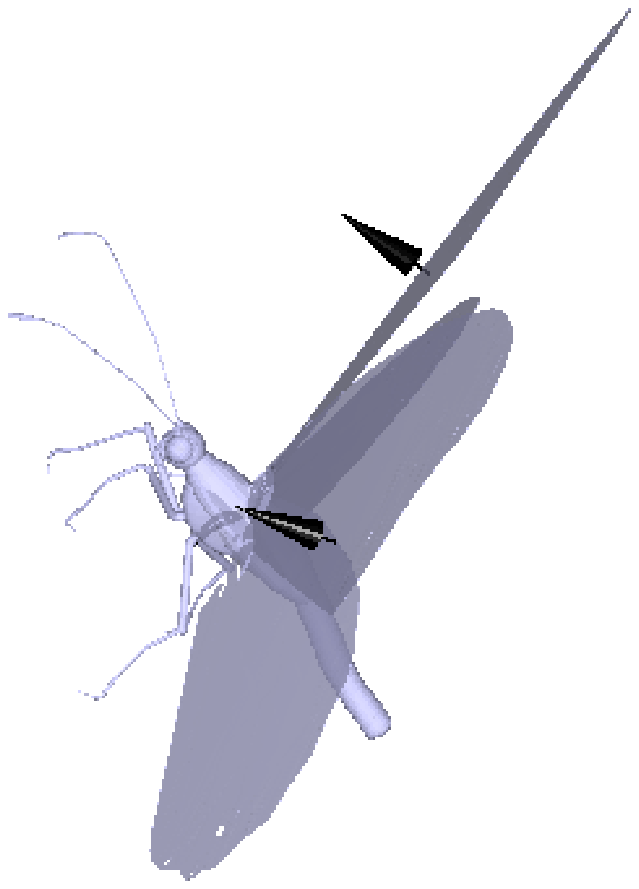}
		}
		\hfill
		\subfigure[$ t=5T/3 $]{
			\includegraphics[trim={13cm 6cm 10cm 2cm},clip,width=0.25\linewidth]{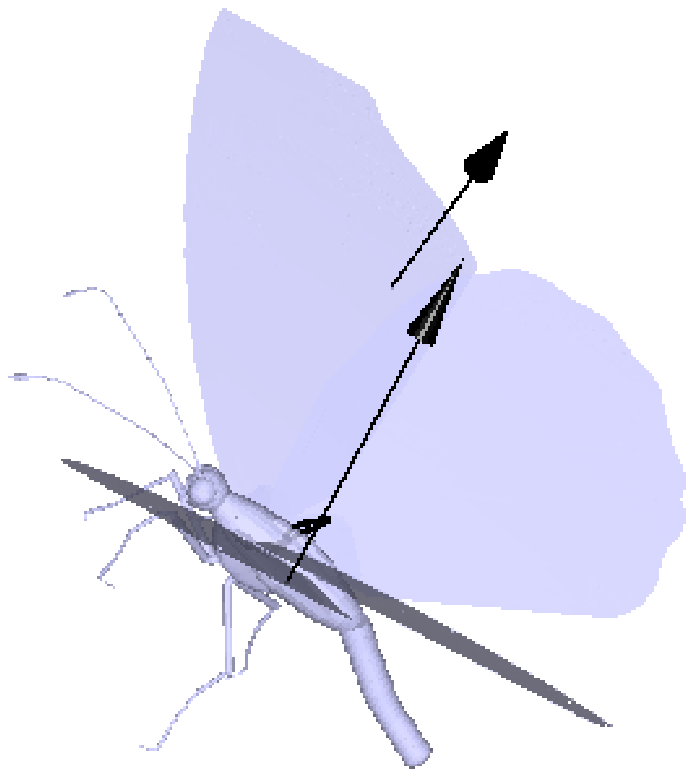}
		}
		\hfill
		\subfigure[$ t=2T $]{
			\includegraphics[trim={11cm 4cm 12cm 4cm},clip,width=0.25\linewidth]{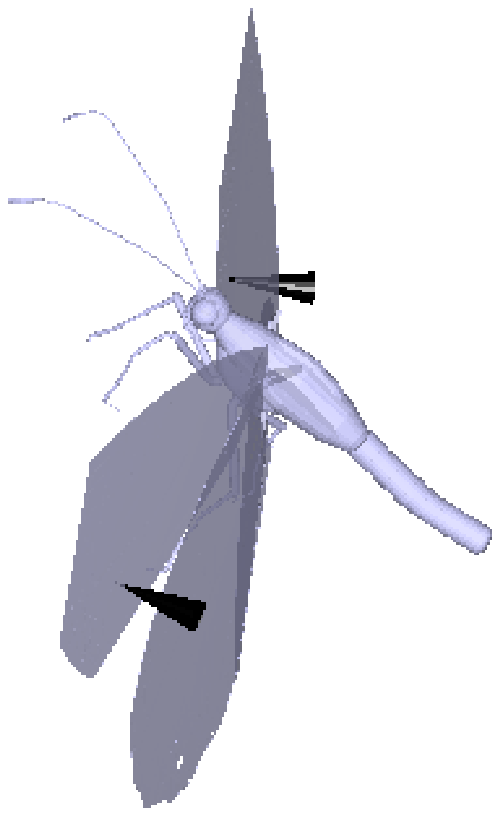}
		}
	}
	\caption{Snapshots of flapping maneuver for hovering: uncontrolled trajectory in (a)-(d); controlled trajectory in (e)-(h) for the same initial condition}
	\label{fig:snapshots}
\end{figure}

The initial states are taken to be,
\begin{gather*}
	x(0) = \begin{bmatrix}
		-0.0003\\
		0.0004\\
		-0.0004
	\end{bmatrix}, \quad
	R(0) = \begin{bmatrix}
		0.7365 & 0.0163 & 0.6763 \\
		-0.0130 & 0.9999 & -0.0100 \\
		-0.6764 & -0.0014 & 0.7366
	\end{bmatrix} \\
	\dot x(0) = \begin{bmatrix}
		-0.2412\\
		0.0100\\
		-0.0787
	\end{bmatrix}, \quad
	\Omega(0) = \begin{bmatrix}
		-0.0437\\
		-2.2907\\
		-0.0487
	\end{bmatrix}.
\end{gather*}
This optimization problem is numerically solved using \textit{fmincon} in MATLAB.
The resulting optimal trajectory errors and the snapshots are illustrated in Figure \ref{fig:hover_control} and \ref{fig:snapshots}, with comparisons to another case without any control.
It is shown that uncontrolled trajectories quickly diverge from the periodic orbit,
whereas the controlled trajectories asymptotically converge to the hovering flight. 

\section{Conclusions}

This paper presents an intrinsic formulation of a Lagrangian system on a Lie group, where the Lagrangian is composed of a configuration-dependent kinetic energy and a potential energy. 
This is utilized for the dynamics of a flapping-wing UAV inspired by Monarch butterfly. 
Two optimization problems are formulated to identify a periodic motion for hovering and also to stabilize it in the framework of model predictive controls. 
Future work includes constructing data-driven feedback control schemes by integrating a set of optimal trajectories computed by the proposed approach for varying initial conditions. 


\end{document}